\documentclass[preprint,review,sort&compress,3p]{elsarticle}
\usepackage[utf8]{inputenc} 
\usepackage{xcolor}
\usepackage{colortbl}
\usepackage{CJKutf8}
\usepackage{mathrsfs}
\usepackage{amsfonts}
\usepackage{cancel}
\usepackage{graphicx}
\usepackage{subcaption}
\usepackage{amssymb}
\usepackage{hyperref}
\usepackage{multirow}
\usepackage{amsthm}
\usepackage{amsmath}
\usepackage{epstopdf}
\usepackage{verbatim}
\usepackage{booktabs}%
\usepackage{flafter}
\usepackage{algorithm}
\usepackage{algorithmicx}
\usepackage{algpseudocode}
\usepackage{makecell}
\usepackage{caption,subcaption}
\usepackage{multirow}
\usepackage{appendix}
\usepackage{float}
\usepackage{array, caption, threeparttable}
\usepackage[font=small,labelfont=bf,labelsep=none]{caption}
\numberwithin{equation}{section}

\journal{ }
\begin{document}
\begin{CJK}{UTF8}{gbsn}
\begin{frontmatter}
\title{MgFNO: Multi-grid Architecture Fourier Neural Operator for Parametric Partial Differential Equations}
\cortext[cor1]{Corresponding author}
\author[a]{Zi-Hao Guo}
\author[a]{Hou-Biao Li\corref{cor1}}
\ead{lihoubiao0189@163.com}

\address[a]{School of Mathematical Sciences, University of Electronic Science and Technology of China, Chengdu, 611731, P. R. China}
\date{May 14, 2025}

\begin{abstract}
 Neural operators are a new type of models that can map between function spaces, allowing trained models to emulate the solution operators of partial differential equations(PDEs). This paper proposes a multigrid Fourier neural operator (MgFNO) that accelerates the training of traditional Fourier neural operators through a novel three-level hierarchical architecture. The key innovation of MgFNO lies in its decoupled training strategy employing three distinct networks at different resolution levels: a coarse-level network first learns low-resolution approximations, an intermediate network refines the solution, and a fine-level network achieves high-resolution accuracy. By combining the frequency principle of deep neural networks with multigrid methodology, MgFNO effectively bridges the complementary learning patterns of neural networks (low-to-high frequency) and multigrid methods (high-to-low frequency error reduction).Experimental results demonstrate that MgFNO achieves relative errors of 0.17\%, 0.28\%, and 0.22\% on the Burgers' equation, Darcy flow, and Navier-Stokes equations, respectively, representing reductions of 89\%, 71\%, and 83\% compared to the conventional FNO. Furthermore, MgFNO supports zero-shot super-resolution prediction, enabling direct application to high-resolution scenarios after training on coarse grids. This study establishes an efficient and high-accuracy new paradigm for solving complex PDEs dominated by high-frequency dynamics. Code and data used are available on \url{https://github.com/guozihao-hub/MgFNO/tree/master}.
\end{abstract}

\end{frontmatter}

\section{Introduction}
Partial differential equation problems (PDEs) arise in many important application fields, including design optimization, uncertainty analysis, and optimization control. We define PDEs as follows:
\begin{equation}\label{1.1}
\left
\{\begin{matrix}
 &\mathcal{A}u=f,~~~~in~  \Omega,\\
 & u=u_{b},~~~~on~\partial\Omega ,
\end{matrix}
\right.
\end{equation}
where $\Omega \in \mathbb{R}^{d}$, $d\in \mathbb{N}_{+}$. $\mathfrak{A}$, $\mathfrak{F}$ are two function spaces that exist on $\Omega$, and $\mathfrak{A}_{b}$ is a function space on $\partial\Omega $. $u\in \mathfrak{A}$, $f\in \mathfrak{F}$ and $\mathcal{A}:\mathfrak{A}—\to \mathfrak{F}$ is a nonlinear differential operator. When the partial differential equation is discretized (finite element, finite difference and finite volume \cite{FiniteElementMethod,FiniteDifferenceMethod,FiniteVolumeMethod}), the above partial differential equation will become a linear system 
\begin{equation}\label{1.2}
    Au=f.
\end{equation}
This linear system is usually large-scale in practice and often solved using iterative methods. Commonly employed iteration methods include splitting iterative methods such as Jacobi iteration, Gaussian-Seidel iteration, and successive over-relaxation (SOR) iteration, as well as Krylov subspace methods such as the conjugate gradient method (CG) and the generalized minimal residual (GMRES) method  \cite{GMRES}.

Traditional numerical methods for PDEs solve equations through discrete space. Consequently, the traditional numerical solution of partial differential equations entails a trade-off in the selection of discrete resolution. That is to say, a coarse grid (low resolution) has a fast calculation speed but low accuracy, whereas a fine grid (high resolution) has high accuracy but slow calculation speed. In practice, for complex partial differential equation systems, very fine discretization is usually required, which requires a huge amount of time consumption and computational resources, which brings practical difficulties to the solution of partial differential equations \cite{2015Convergence}.

The multigrid method (MG) \cite {MG1, MG2} is one of the most remarkable classical methods. The computational complexity of the MG method for Poisson's equation is only $O(n)$, where $n$ is the size of the matrix. Another popular high-performance numerical method is the spectral method, which has a complexity of $O(nlogn)$.

Recently, with the explosive development of computing resources and data sampling technology, machine learning and data-driven methods, especially deep neural networks, have rapidly developed in the field of information data and have been widely applied to image recognition and natural language processing problems \cite{ImageRecognition, nlp}. However, when analyzing complex physical, biological, or engineering systems, the cost of data collection is often prohibitively high, and we inevitably face the challenge of making decisions based on partial information. In small data systems, the vast majority of state-of-the-art machine learning techniques (such as deep/convolutional/recursive neural networks) lack robustness and cannot provide any guarantee of convergence. To address this issue, Physical Information Neural Networks (PINNS) \cite{PINNS} is proposed, which directly parameterizes the solution function into a neural network. Palha and Rebelo 
 \cite{PINNS} point out that for many physics related modeling problems, there is often a large amount of prior knowledge that has not yet been used in modern machine learning practices. These prior knowledge can play the role of regularization constraints, constraining the space of tolerable solutions to a manageable range. As a reward, encoding this structured information into the learning algorithm amplifies the data seen by the algorithm, enabling it to quickly converge to the correct solution, and even in situations where only a small amount of training data can be obtained, the machine algorithm still maintains good generalization ability \cite{PINNS2}.

Machine learning methods can provide fast solvers that approximate or enhance traditional disciplines, thereby becoming the key to the scientific discipline revolution. However, PINNS maps between finite dimensional spaces.  For a new instance of a differential system, this method requires retraining the model and is limited to the setting of known potential PDEs. In response to the above issues, Lu and pang \cite{DeepONet} proposed the Operator neural network (DeepONet), which compensates for the mesh dependency of the finite dimensional operator method by generating a single network parameter set that can be used with different discretizations. It can transfer solutions between grids of different resolutions. In addition, operator neural networks only need to be trained once. The solution for obtaining new instances of parameters only requires forward propagation of the network, thereby alleviating the main computational problems that arise in neural finite element methods. Finally, operator neural networks do not require knowledge of underlying PDEs, only data, and are completely data-driven methods. Due to the high cost of calculating integral operators, DeepONet has not yet produced effective numerical algorithms, and therefore cannot achieve successful replication of convolutional or recursive neural networks in finite-dimensional environments. In a different approach the graph neural operator (GNO) was proposed \cite{GNO}. The GNO focuses on learning of the infinite-dimensional mapping by composition of nonlinear activation functions and certain class of integral operators. Subsequently, Li, Kovachki and Azizzadenesheli 
 \cite {FNO1} proposed the fourier neural operator (FNO) based on this foundation. FNO uses fast fourier transform (FFT) to approximate partial differential equations, which is more in line with the laws of physical equations than convolutional and recursive neural networks, and has shown good performance and speed in flow problems such as Darcy flow and N-S equations. Zhang and Zhao \cite{FNO2} has developed a deep learning based FNO model to solve three types of problems in PDE control of two-dimensional oil/water two-phase flow underground. His work has shown great potential in replacing traditional numerical methods with neural networks for reservoir numerical simulation. 

\begin{figure}[h]
    \centering
    \includegraphics[width=1.0\textwidth]{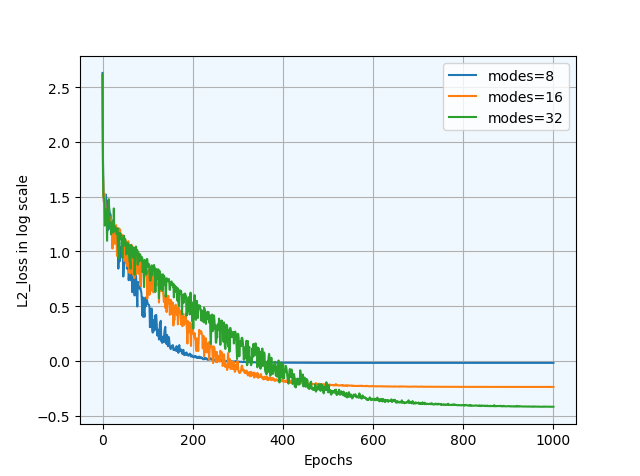}
    \vspace{0.5em}
    \caption{ :  FNO introduces a trade-off between the speed of computation and the accuracy of fitting with different "modes" according to the parameter settings in \cite{FNO1} to solve 1-d Burger's problem.}
    \label{fno_different_modes}
\end{figure}

When constructing a network model, FNO needs to specify the parameter "modes", which represents the number of frequency components retained during the Fourier transform. If only a small number of low-frequency components are retained, the training time of the neural network model is short, but the training error is often large. If more frequency components are retained, the training error is relatively small, but the model training time is often longer. Figure \ref{fno_different_modes} illustrates the trade-off between computation speed and fitting accuracy imposed by FNO. Another shortcoming of FNO is that the basis functions in the Fourier transform are generally frequency localized. This means that a specific basis function corresponds to a specific frequency, resulting in clear and concentrated basis functions in the frequency domain. While the Fourier transform is effective at analyzing the frequency components of signals, it is not suitable for representing local transformations of signals in time or space. This limitation has led to the development of other transformations, such as the wavelet transform. The wavelet neural operator (WNO) \cite{WNO} is proposed to overcome this aforementioned shortcoming. It transforms the snapshots into their high and low frequency components using discrete wavelet transform.
\hspace*{\fill}

A dynamic system viewpoint has been explored in several papers, such as \cite{FPrinciple1, FPrinciple2, FPrinciple3}, to understand the frequency tendency of the neural network training process. Experiments conducted on synthetic low-dimensional data demonstrate that during training, deep neural networks (DNNs) often fit target functions from low to high frequencies. This inherent frequency bias is referred to as the frequency principle (F-Principle) \cite{FPrinciple1}. While FNO mitigates the harm caused by the F-Principle by controlling the number of frequency "modes," it does result in some loss of accuracy during model training.

In this paper, our goal is to enhance the training process of FNO by focusing on frequency aspects.
To begin with, we propose a three-layer V-cycle multigrid architecture for training the model. This involves training the model multiple times on a coarse grid and then transferring it to a fine grid to expedite the training process. The use of multigrid effectively addresses the challenge of eliminating high-frequency components. The underlying concept is that high-frequency components on fine grids can be transformed into low-frequency components on coarse grids. Existing research has already explored the integration of multigrid methods to develop efficient DNNs. For example, in \cite{MgNet}, a unified model called MgNet is introduced, which combines convolutional neural networks (CNN) for image classification with multigrid (MG) methods for solving discretized partial differential equations (PDEs). Another study \cite{MgNet2} proposes a linear U-net structure as a solver for linear PDEs on a regular mesh. Additionally, we incorporate compact support activation functions to ensure that FNO achieves rapid and consistent convergence across multiple scales \cite{MscaleDNN}.
Using the aforementioned observations and concepts, we have developed a novel Fourier neural operator known as MgFNO. This innovative approach addresses the issue of slow training speed in FNO from a frequency perspective. Moreover, our proposed MgFNO has proven to outperform the traditional Fourier neural operator in solving various equations, including the 1-D Burgers' equation, 2-D Darcy flow, and 2-D time-dependent Navier-Stokes equation.

The remaining sections are organized as follows. In section 2, we introduce necessary notations for the rest of the paper and review the use of the MG method and FNO to solve parametric PDEs. In section 3, we discuss the frequency principle and various methods to overcome low-frequency trends. We pay special attention to the construction of a multiscale DNN that can handle the challenge of converging high-frequency components for neural networks. In section 4, we introduce the proposed MgFNO. Section 5 presents the numerical experiments and comparisons among the classical FNO methods. Finally, concluding remarks are given in section 6.

\section{The MG method and FNO to solve parametric PDEs}
In this section, we briefly describe multigrid method and fourier neural operator to solve the PDEs. MG method is one of the most efficiency methods for solving PDEs. FNO is a new and effective framework for solving PDEs, which surpasses the existing neural network methods in the benchmark PDE problems. 

\subsection{The Multigrid method}
Before using the MG method to solve PDEs (\ref{1.1}), it is necessary to discretize the PDEs. The finite difference method (FDM)is one of the most commonly used discretization methods, which directly approximates differential equations on spatial grids. By using differential discretization, each node is approximated and the differential operator is replaced by differential approximation. In the finite element method (FEM), the solution area is first divided into a finite number of small elements, and appropriate mathematical interpolation functions are used within each element to approximate the solution. By introducing shape functions and element stiffness matrices, the behavior of the overall system can be transformed into relationships between local elements.  The finite volume method (FVM) focuses on conserving fluxes by integrating the governing equations over discrete control volumes. This approach ensures that the flux entering and leaving a control volume is accurately accounted for, preserving the overall balance in the system. In summary, both of these methods discretize a continuous problem into a discrete problem, and then obtain the numerical solution of the problem through numerical methods.

For linear system (\ref{1.2}), iterative method is one of the basic numerical methods. Let the initial guess value of the solution be $u_{0}$, and update scheme is represented by $\mathscr{H}$. We can write the iterative method as
\begin{equation}
    u_{t+1}=u_{t}+\mathscr{H}(f-Au_{t}), t=0,1,...,T.
    \label{Iterative method}
\end{equation}
 Jacobi iteration is a classic iterative method, which has the form
\begin{equation}
    u_{t+1}=u_{t}+D^{-1}(f-Au_{t}), t=0,1,...,T.
    \label{jacobi}
\end{equation}
where $D$ is the diagonal matrix of $A$. We utilize the Jacobi iteration  to solve the 1-D Poisson equation
\begin{equation}
    \nabla ^2u=f.
    \label{possion}
\end{equation}
and utilize central difference for discretization. Eq. (\ref{1.2}) can be written in the pointwise form
\begin{equation}
    u_{i}^{(k+1)} = u_{i}^{(k)} + \frac{1}{2} \left( u_{i-1}^{(k)} + u_{i+1}^{(k)} - 2u_{i}^{(k)} -  f_i \right).
    \label{pointwise_form}
\end{equation}
Let $u_i$ be the true solution in point $i$, and define error $e_i^{k+1}=u_i^{k+1}-u_i$. Then we have
\begin{equation}
    e_i^{(k+1)} = e_i^{(k)} + \frac{1}{2} \left( e_{i-1}^{(k)} + e_{i+1}^{(k)} - 2e_i^{(k)} \right).
    \label{pointwise_form2}
\end{equation}
The error's fourier transform is $e_{i}^k = \sum_{p=-\infty}^{\infty} c^{k} e^{j2\pi px_i}$ and $|c^{k}|$ is the amplitude of the corresponding frequency. Then we define the convergence factor of Jacobi iteration in Poisson equation is $ \mu_{loc} =|\frac{c^{k+1}}{c^{k}}|$.
\begin{figure}[H]
    \centering
    \includegraphics[width=0.7\textwidth]{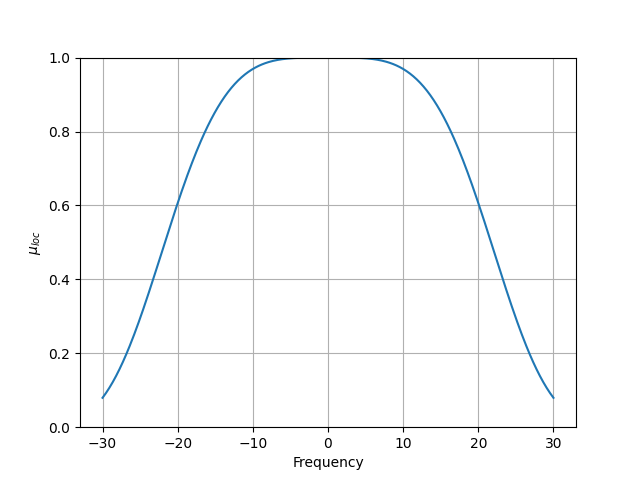}
    \caption{ :  Distribution of the convergence factor $\mu_{loc}$.}
    \label{jacobi}
\end{figure}

Figure \ref{jacobi} illustrates the distribution of convergence factors for the Jacobi iterative method when solving linear systems derived from the Poisson equation. It can be clearly observed that the iterative method (\ref{Iterative method}) effectively reduces high-frequency errors but generally performs poorly in handling low-frequency errors. This limitation arises from the localized nature of iterative methods, which struggle to eliminate low-frequency global error components generated by overall patterns. Notably, the frequency characteristics of iterative methods exhibit opposite convergence tendencies compared to Fourier neural operators. Therefore, we aim to explore improvements to iterative methods, thereby enhancing the performance of Fourier neural operators.\\
To address the limitations of iterative methods (\ref{Iterative method}), the multigrid (MG) method was initially proposed by R. P. Fedorenko in the 1960s and subsequently refined by Achi Brandt and others in the 1970s. The fundamental principle of this approach lies in utilizing hierarchical grid structures to accelerate convergence through two key algorithmic steps: (1) fine grid smoothing and (2) coarse grid correction. In the fine grid smoothing phase, conventional iterative methods (smoothers) are employed to effectively eliminate high-frequency approximation errors. However, recognizing that such local operations remain ineffective for low-frequency error components, the method strategically employs coarse grid correction - where low-frequency errors are approximated and corrected through solutions computed on coarser grids. This dual-scale mechanism enables the MG method to simultaneously achieve both rapid error reduction through coarse grid processing and high precision through fine grid refinement, with the complete algorithmic procedure detailed in Algorithm \ref{alg1}.

\begin{algorithm}[H]
        \caption{MG method\cite{MG1}}
        \label{alg1}
        \begin{algorithmic}[1]    
        \Require coefficient matrix $A$, right hand item $f$, initial guess value of the solution $u^{1,0}$, number of grid layers $J$, number of smoothness on each grid layer $v$  
        \Ensure approximate solution $u$    
        \State Initialization：$$f^{1}=f,A^{1}=A,u^{1,0}=0$$
        \State Restriction from fine to coarse level
        \For{each $l \in [1,J]$}
            \If{$l=1$}
                \For{each $i \in[1,v]$}
                    \State $u^{1,i}=Smoothing(A^{1},u^{1,i-1},f^{1})$
                \EndFor
                \State $r^{1}=f^{1}-A^{1}u^{1,v},A^{1}e^{1}=r^{1}$
            \EndIf
                            
            \If {$l\in[2,J-1]$}
                \State $A^{l}=Restriction(A^{l-1}),r^{l}=Restriction(r^{l-1}),e^{l,0}=0$
            
                \For{each $i \in[1,v]$}
                    \State $e^{l,i}=Smoothing(A^{l},e^{l,i-1},r^{l})$
                \EndFor
                \State $r^{l}=r^{l}-A^{l}e^{l,v}$
            \EndIf
            \If {$l=J$}
                \State $A^{J}=Restriction(A^{J-1}),r^{J}=Restriction(r^{J-1}),e^{J,0}=0$
                \State $e^{J,*}=(A^{J})^{-1} r^{J}$
            \EndIf             
        \EndFor
        \State  Prolongation from coarse to fine level
        \For{$l\in[J-1,1]$}
            \State{$e^{l,*}=e^{l,v}+ Prolongation(e^{l+1,*})$}
        \EndFor\\
        \Return approximate solution $u$
    \end{algorithmic}
\end{algorithm}

\subsection{Neural Operator }
Neural networks can not only approximate any continuous function but also any nonlinear continuous operator \cite{Rossi_Conan-Guez_2005} 
 (a mapping from a space of functions into another space of functions). Before reviewing the theorem of operators, we introduce some notation. Let $D \in \mathbb{R}^d$ be the $n$-dimensional domain with boundary $\partial D$. For a fixed domain $ D=(a,b)$ and $x\in D$, consider a PDE which maps these function spaces containing the source term $f(x,t):D\to \mathbb{R}$, the initial condition $u(x,0):D\to \mathbb{R}$, and the boundary conditions $u(\partial D,t):D\to \mathbb{R}$ to the solution space $u(x,t):D\to \mathbb{R}$, with $t$ being the time coordinate. Then let us define two complete normed vector spaces ($\mathcal{A},||\ldotp ||_{\mathcal{A}}$) and ($\mathcal{U},||\ldotp ||_{\mathcal{U}}$) of functions taking values in $\mathbb{R}^{d_a}$ and $\mathbb{R}^{d_u}$. These function spaces are denoted as $\mathcal{A}(D;\mathbb{R}^{d_a})$ and $\mathcal{U}(D;\mathbb{R}^{d_u})$, which are called as Banach spaces.

Let $a\in\mathcal{A}$, $u\in\mathcal{U}$ be two input and output functions existing respectively in the function spaces mentioned earlier. Further, consider that $u_i=\mathcal{G}(a_i)$ and we have access to $N$ number of pairs $(a_j,u_j)_{j=1}^N$, where $\mathcal{G}:\mathcal{A}\to \mathcal{U}$ be a non-linear map between the two function spaces. We aim to learn this mapping using a neural network $\mathcal{G}^{+}(a)=u,\mathcal{G}^{+}\approx \mathcal{G}$ from the pairs $(a_j,u_j)_{j=1}^N$.

 The MG method solves the solution $u(x,t)$ directly which satisfies the differential operator $\mathcal{G}$ in the domain $D$ and boundary conditions $\partial D$.
 So when the parameters of $a\in\mathcal{A}$ change, the MG method must recalculate the solution of the PDEs. In contrast, the neural operator directly learn the nonlinear differential operator $\mathcal{G}$ so that the solutions to a family of PDEs can be obtained for different sets of input parameters $a\in\mathcal{A}$.

\subsection{Fourier Neural Operator}
The network layer of neural operator is as follows:
\begin{equation}
    {v_0(x)=P(a(x)),\quad v_{t+1}=\sigma (Wv_{t}(x)+(\mathcal{K}(a;\phi )* v_t)(x)),\quad  u(x)=Qv_T(x) }
    \label{fno_eq}
\end{equation}
\begin{equation}
    (\mathcal{K}(a;\phi )* v_t)(x)=\int_Dk(a(x,y),x,y;\phi)v_t(y)dy,\forall x\in D.
    \label{fno_eq2}
\end{equation}

Firstly, $P:\mathbb{R}^{d_a}\to\mathbb{R}^{d_v}$ lifts the inputs $a(x,t)\in \mathcal{A}$ to  the high dimensional space $v_0$ ,which takes values in $\mathbb{R}^{d_v}$. The local transformation $P$ can be achieved by constructing a fully connected neural network. Then,  perform $T$ iterations on $v_0 (x)$. The step-wise updates are defined as Eq.(\ref{fno_eq}).
Where $W:\mathbb{R}^{d_v}\to\mathbb{R}^{d_v}$ is a linear transformation, $\mathcal{K}(a;\phi )$  is the integral operator  and $\sigma$ is the non-linear activation function. Finally, when all the updates are performed, we define another local transformation $Q:\mathbb{R}^{d_v}\to\mathbb{R}^{d_u}$ to transform back the output $v_T(x)$ to the solution $u(x)$.

According to schwartz kernel theorem \cite{schwartzKernel}, for every such continuous kernel integral operator $\mathcal{K}(a;\phi)$ parameterized by $\phi \in \theta$ there exists one and only one distribution $k$. The kernel integral operator $\mathcal{K}(a;\phi)$ is defined as Eq.(\ref{fno_eq2}), where $k_{\phi}:\mathbb{R}^{2d+d_a}\to\mathbb{R}^{d_v\times d_v }$ is a neural network parameterized by $\phi \in \theta$. Since the kernel $k$ and $\mathcal{K}$ define the infinite-dimensional spaces, Eq.(\ref{fno_eq2}) can learn the mapping of any infinite-dimensional function space. There are various types of kernel neural networks, including fully connected neural networks, graph convolutional networks, and so on.

\begin{figure}[h]
    \centering
     \vspace{0.5cm}
    \includegraphics[width=1.0\textwidth]{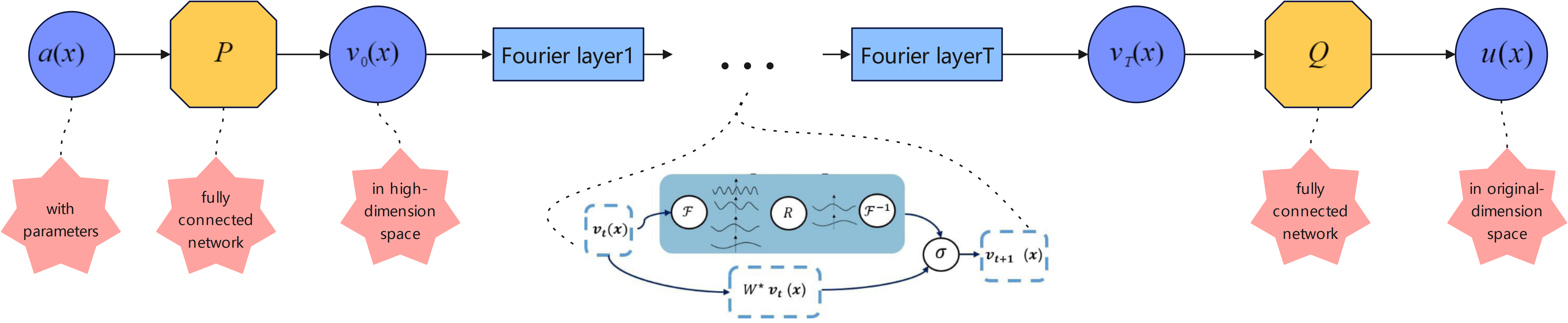}
    \vspace{0.2cm}
    \caption{： FNO network: The input PDE parameter $a$ is passed from $P$ to $v_0(x)$, then it undergoes iteration in the Fourier space to obtain $V_T(x)$, and finally it goes through $Q$ to reach the spatial dimension of the solution $u$.}
    \label{Self drawn fno}
\end{figure}

Let $\mathcal{F}$ is the fourier transform, $\mathcal{F}^{-1}$ is the inverse fourier transfrom. According to the fourier convolution theorem \cite{FourierConvolutionTheorem}, the integral function $\int_Dk(\phi)v_t(y)dy$ is converted to the product form 
$\mathcal{F}^{-1}(\mathcal{F}(k_{\phi})\cdotp\mathcal{F}(v_t) )(x)
$.\\$\mathcal{F}^{-1}(\mathcal{F}(k_{\phi})\cdotp\mathcal{F}(v_t) )(x)$ is the key to learn from parameter PDE in fourier space. The high fourier modes are removed to reduce neural network training time, and the filtering strategy $R$ is to preserve low-order modes in traditional FNO. The traditional FNO algorithm believes that most of the energy of the fourier transform is mainly concentrated in the four corners of the coefficient matrix, namely the low order region. Therefore, removing high-frequency components can improve the model's generalization ability and speed. The framework of FNO is presented in Figure \ref{Self drawn fno}.
Let the parameters of network $\mathcal{G}^{+}$ are $\theta_{NN}$, for training this network, FNO defines an appropriate loss function as follows:
\begin{equation}
    \theta _{NN}=argminLoss(\mathcal{G}^{+}(a),u).
    \label{fno_eq3}
\end{equation}
We try to learn the parameter space (weights-$W$, biases-$b$ and the filtering strategy $R$) by approaching to the loss function from a data-driven perspective.

\section{Frequency Tendency in Deep Learning and Compensation Methods}
Although deep neural networks (DNNs) possess powerful capabilities, their theoretical analysis is not yet sufficient, so DNNs are also known as "black boxes". However, in recent years, some studies have revealed this magical "black box" by demonstrating the frequency principle of deep neural network training behavior. 
In this section, we first introduced the low-frequency bias in the neural network fitting process. Subsequently, we will explore the utilization of multiscale neural networks and the conversion of coarse and fine grids as approaches to overcome the high-frequency curse.

\subsection{Frequency principle}
\begin{figure}[h]
    \centering
    \begin{minipage}[b]{0.3\textwidth}
        \centering
        \includegraphics[width=\textwidth]{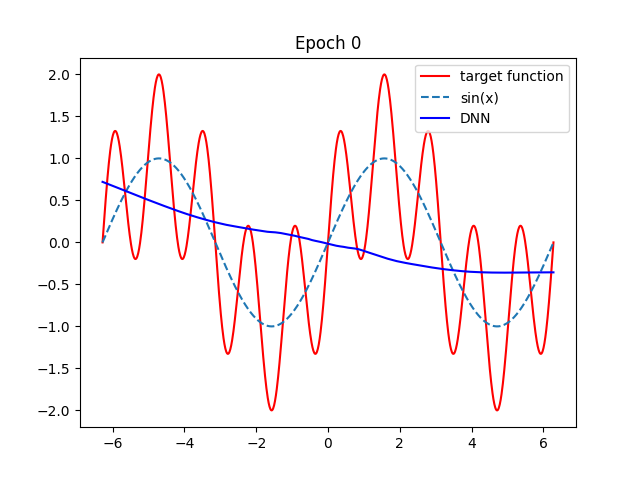}
        
    \end{minipage}
    \hfill
    \begin{minipage}[b]{0.3\textwidth}
        \centering
        \includegraphics[width=\textwidth]{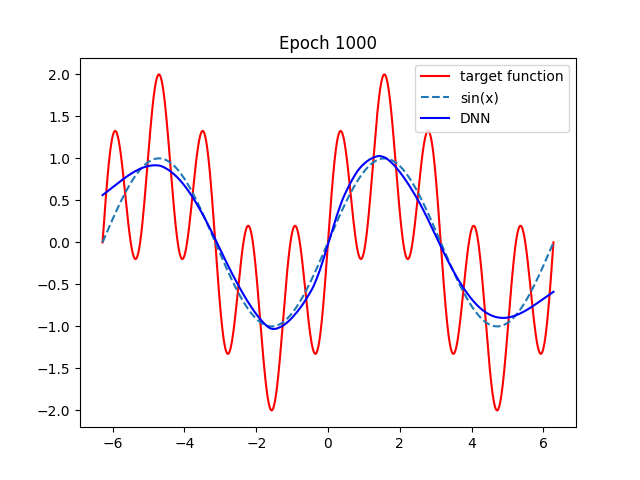}
        
    \end{minipage}
    \hfill
    \begin{minipage}[b]{0.3\textwidth}
        \centering
        \includegraphics[width=\textwidth]{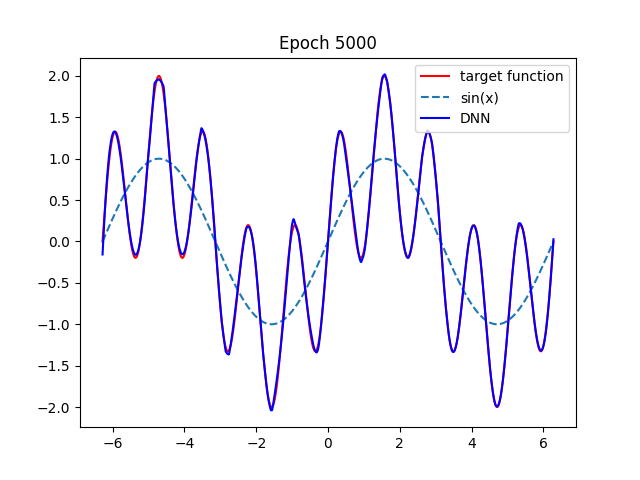}
        
    \end{minipage}
    \caption{: Training process of DNNs in the spatial domain.}
    \label{Fprinciple 1}
\end{figure}

\begin{figure}[h]
    \centering
    \subcaptionbox{\label{fig:image1}}%
        [.4\linewidth]{\includegraphics[width=\linewidth]{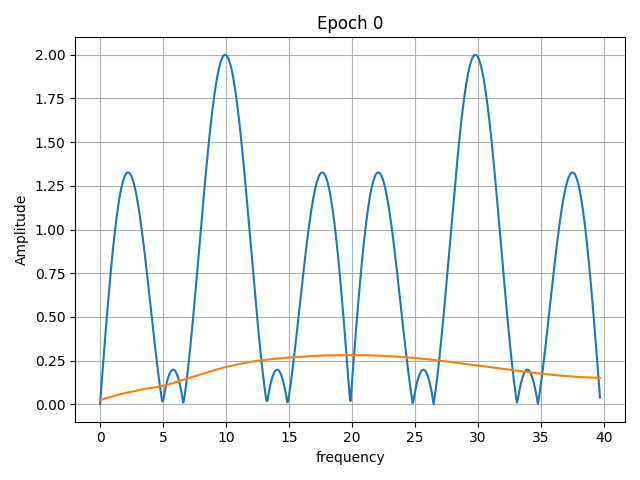}}\hfill
    \subcaptionbox{\label{fig:image2}}%
        [.4\linewidth]{\includegraphics[width=\linewidth]{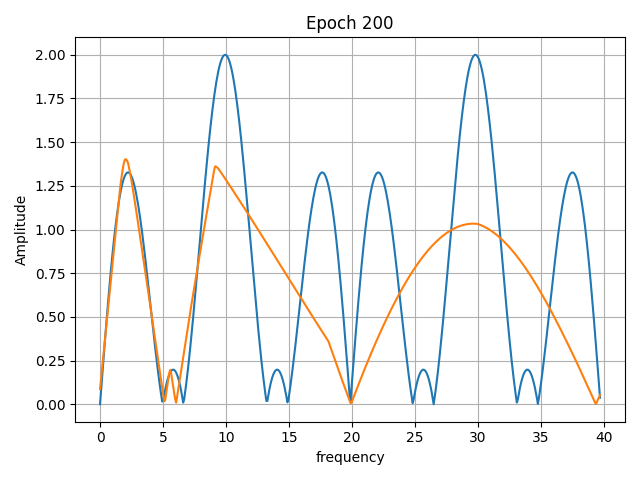}}
    
    \vspace{0.2cm} 
    
    \subcaptionbox{\label{fig:image3}}%
        [.4\linewidth]{\includegraphics[width=\linewidth]{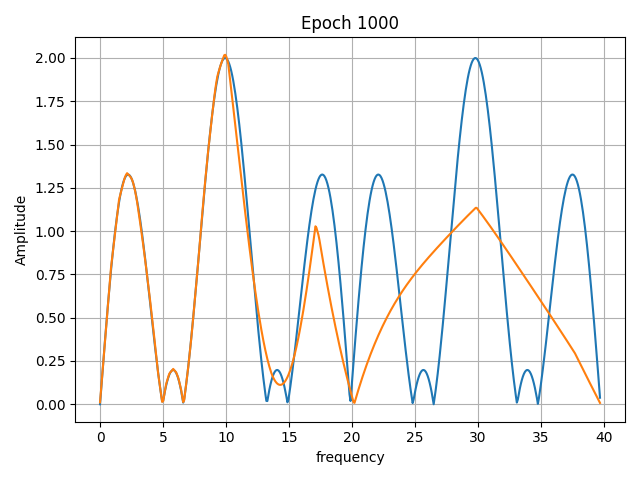}}\hfill
    \subcaptionbox{\label{fig:image4}}%
        [.4\linewidth]{\includegraphics[width=\linewidth]{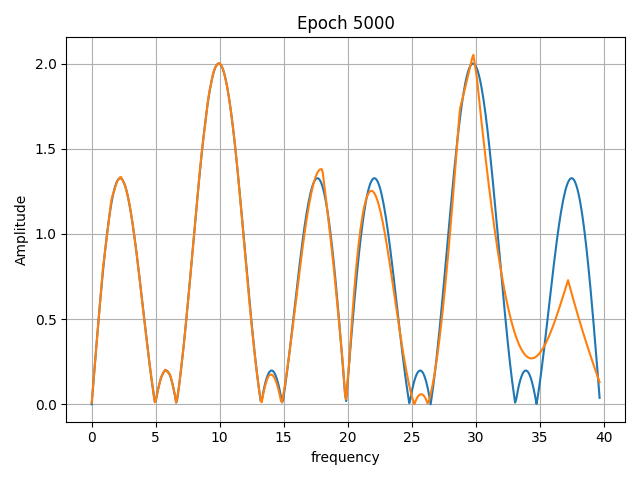}}

    \caption{: Training process of DNNs in the frequency domain. Amplitude $|\hat{y}_k|$ vs. frequency.}
    \label{Fprinciple 2}
\end{figure}

Recently, several papers\cite{FPrinciple1,FPrinciple2,FPrinciple3} have discovered that DNNs exhibit frequency-dependent convergence behavior. Specifically, during training, DNNs tend to fit target functions from low to high frequencies. This phenomenon is known as the F-Principle of DNNs\cite{FPrinciple1}. The slow convergence of DNNs in the later stages of training may be attributed to the challenges in learning high-frequency components of the data. In order to demonstrate the F-Principle, we utilize 1-D synthetic data to observe the evolution of a DNN in both the frequency and spatial domains. Our objective is to train a DNN to fit a 1-D target function $f(x) = sinx + sin5x$, which comprises two frequency components. The data set consists of 1,000 evenly distributed sampling points within the interval $[-2\pi,2\pi]$. We compute the discrete Fourier transform of $f(x)$ using the formula $\hat{f}_k=\frac{1}{n}\textstyle\sum_{i=0}^{n-1}f(x_i)e^{-j2\pi ik/n}$, where $k$ represents the frequency. We initialize the parameters of the DNN using a Gaussian distribution with a mean of 0 and a standard deviation of 0.1. We employ a ReLU-DNN with widths of 1-200-200-1, and set the learning rate to 0.0005. Figure \ref{Fprinciple 1} and Figure \ref{Fprinciple 2} illustrate the rapid convergence of the DNN for low frequencies, while the convergence is significantly slower for high frequencies.

It is worth noting that the iterative method (\ref{Iterative method}) to solve linear systems often converges first on the high-frequency component of the error, and then on the low-frequency component. In other words, the convergence tendency of iterative methods is exactly opposite to that of DNNs. The MG method improves upon the iterative method by utilizing the characteristic that the low-frequency components on the coarse grid are the high-frequency components on the fine grid. This has inspired us, and we have also attempted to use this feature to construct a multigrid model to accelerate the fitting process of FNO.

Next, we will define a neural network and explain the concept of the frequency principle. Firstly, let $T(x)=Wx+b$ be a linear function, where $W=(w_{ij})\in\mathbb{R}^{m\times n}$ represents the weights and $b\in\mathbb{R}^{m}$ represents the biases. Then, let $\sigma(u):\mathbb{R}\rightarrow \mathbb{R}$ be a nonlinear activation function. Common activation functions include the sigmoid function, Tanh function, ReLU function, etc. Finally, a DNN $NN(x)$ with $L$ layers can be expressed as:
\begin{equation}
    NN(x)=T^{L-1}\odot \sigma \odot T^{L-2}\ldots\odot  T^{1}\odot \sigma \odot  T^{0}(x).
    \label{FNN}
\end{equation}

This DNN (\ref{FNN}) is also referred to as having $L-1$ hidden layers, with the $l$-th layer containing $n^{l}$ neurons. Assuming the true function is $f(x)$, we utilize a neural network $NN(x)$ (\ref{FNN}) to fit it. For simplicity, we denote the parameter vector $\theta$ as all the parameters in the DNN. Assuming we use N data points to train the network, i.e., the training set is $\left \{x_1, x_2, \ldots, x_N \right \}$, the loss function is defined as
$$Loss(\theta )=\displaystyle\sum_{i=1}^{N}|f(x_{i})-NN(x_{i},\theta )|^{2}.$$

We define the Fourier transform of $f(x)$ as $\hat{f}(k)$, where $\hat{f}(k) = \frac{1}{\sqrt{2\pi}}\int_{-\infty}^{+\infty}f(x)e^{-ikx}dx$. Parseval's equality states that energy is conserved in both the time-domain and frequency-domain spaces, i.e., $||\hat{f}(k)|| = \sqrt{2\pi}||f(x)||$. Here, the norm is defined as a 2-norm. According to Parseval's equality, we have
\begin{equation}
    Loss(\theta )=\int_{-\infty }^{+\infty }|f(x)-NN(x,\theta )|^2dx=\int_{-\infty }^{+\infty }|\hat{f}(k)-\hat{NN}(k,\theta )|^2dk.
\end{equation}

We provide a specific theorem for the frequency principle which is proven in \cite{Luo_2021}.

THEOREM 1 (\cite{Luo_2021}). Assuming we use gradient descent to train the aforementioned neural network $NN(x,\theta)$, and $t$ represents the number of training iterations for the neural network. If during the training process, the following assumptions hold:\\
1.\; $\theta (t)\ne$ constant and $sup_{t\geqslant 0}|\theta (t)|\le R$ for a constant $R$;\\
2.\; $inf_{t>0}|\nabla_{\theta }Loss(\theta ) |>0;$\\
Then, in a fixed training time interval $t$, there exists a constant $C$ such that\\
$$\frac{dLoss^{+}_{\eta }(\theta )/dt}{dLoss(\theta )/dt}<C\eta^{-m}, \; and\; \; \frac{dLoss^{-}_{\eta }(\theta )/dt}{dLoss(\theta )/dt}>1-C\eta^{-m}  ~~for~  1\leqq m\leqq2r-1.$$
Wherein, $$Loss^{-}_{\eta }(\theta )=\int_{B_{\eta}}^{}|\hat{f}(k)-\hat{NN}(k,\theta )|dk,\; Loss^{+}_{\eta }(\theta )=\int_{B_{\eta}^{c}}^{}|\hat{f}(k)-\hat{NN}(k,\theta )|dk，$$
where $B_{\eta}$ is the ball with radius $\eta$, $B_{\eta}^{c}$ is its complement.

This theorem reveals the relative rate of change between the low-frequency and high-frequency components of the loss function during training. This theorem indicates that when the gradient descent method is applied to the loss function, the low-frequency part of the loss function converges faster than the high-frequency part.

\subsection{Compensation Methods of Low-Frequency Tendency}
 In order to address the issues caused by slow convergence at high frequencies and speed up the training of FNO, several methods have been proposed. Jagtap and Ameya D \cite{RemedialMethod1} used activation functions $\sigma (\mu ax)$ instead of $\sigma (\mu x)$. Here, $\mu$ is a fixed scaling factor, and $a$ is a trainable parameter for all neurons. Cai and Wei \cite{RemedialMethod2} introduced PhaseDNN, which converts the high-frequency components of data to the low-frequency spectrum for learning, and then converts the learned frequency components back to the original high-frequency spectrum. Subsequently, the multiscale DNN method (MscaleDNN) was proposed \cite{MscaleDNN,RemedialMethod3}, which can be understood as Fourier expansion only in the radial direction. The frequency space conversion can be achieved through scaling, which is equivalent to inverse scaling in spatial space. In this article, we utilize the conversion of coarse and fine grids and MscaleDNN to address the slow convergence issue of high-frequency errors in FNO. Next, let us review the ideas and structure of MscaleDNN in regular DNNs.
 
\begin{figure}[h]
    \centering
    \includegraphics[width=0.7\textwidth]{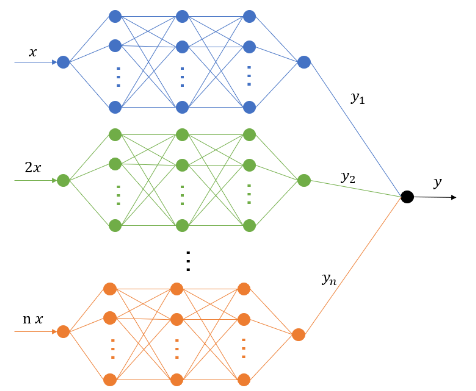}
    \caption{: A common MscaleDNN structure}
    \label{Msnet}
\end{figure}

In MscaleDNN, assuming the function we want to learn is $f(x)$, its fourier transform is $\hat{f}(k),|k|\le K_{max}$.
Divide the frequency of $\hat{f}(k)$ into $M$ segments, that is, the entire frequency domain is divided into the union of $M$ concentric rings $A_i$. 
\begin{equation}
    A_i=\left \{ (i-1)K_0\leq |k|\leq i K_0\right \},K_0=K_{max}/M.
\end{equation}
 So the frequency domain function $\hat{f}(k)$ can be decomposed into the sum of $M$ sub-functions, i.e., $\hat{f}(k)=\displaystyle\sum_{i=1}^{M}\hat{f}_{i}(k)$, and $\hat{f}_{i}(k)$ is a sub-function of $\hat{f}(k)$ corresponding to the frequency band. So corresponding decomposition can also be done in physical space, that is, $f(x)=\displaystyle\sum_{i=1}^{M}f_{i}(x)$, where $f_{i}(x)=\mathcal{F}^{-1}[\hat{f}_{i}(k)](x)$. Then we define $\alpha_i$ as the scaling factor, in order to convert the high frequency region $A_i$ to a low frequency region, i.e., $\hat{f}_i^{scale}(k)=\hat{f}_i(\alpha_ik)$. Similarly， for $f(x)$ in the physical space, 
\begin{equation}
    f_i(x)=f_i^{(scale)}(\alpha_ix).
\end{equation}
A common MscaleDNN structure is shown in Figure \ref{Msnet}. We use neural networks successively in physical space to learn the stretched sub-functions $f_i^{(scale)}(\alpha_ix)$. Furthermore, it is mentioned in \cite{Liu_2020} that when constructing MscaleDNN, activation functions with compact support should be selected. In MgFNO, we select the following function $\phi(x)$ as our activation function.
  \begin{equation}
      \phi(x) =(ReLU(x))^2-3(ReLU(x-1))^2+3(ReLU(x-2))^2-(ReLU(x-3))^2. 
      \label{phi(x)}
  \end{equation}
So $f(x)$ can be expressed as the sum of $M$ sub networks, i.e., $f(x)\approx \displaystyle\sum_{i=1}^{M}NN_{i}(\alpha_i x)$.  

\hspace*{\fill}

\begin{figure}[H]
    \centering
    \begin{minipage}[b]{0.90\textwidth}
        \centering
        \includegraphics[width=\textwidth]{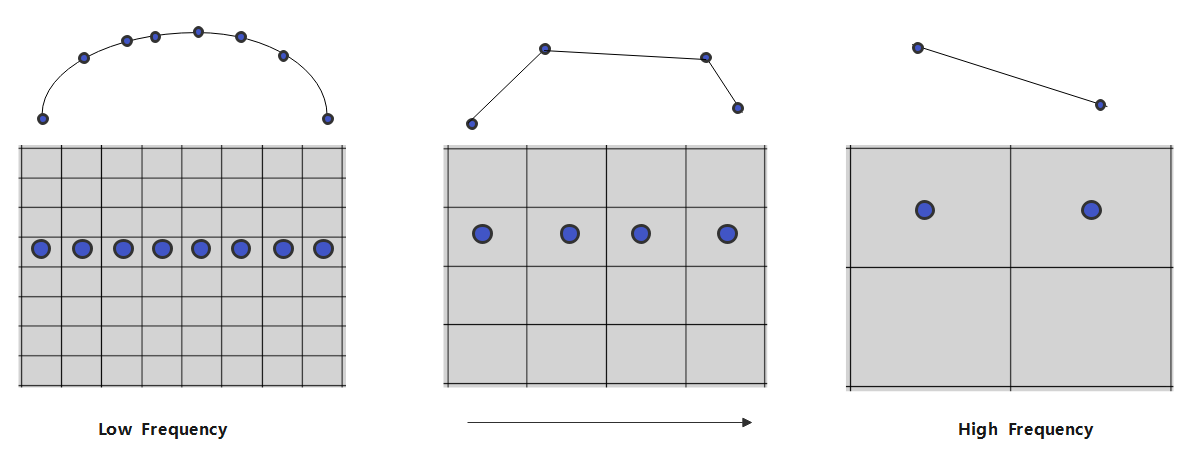}     
    \end{minipage}
    \vspace{0.8cm}
  \begin{minipage}[b]{0.90\textwidth}
        \centering
       \includegraphics[width=\textwidth]{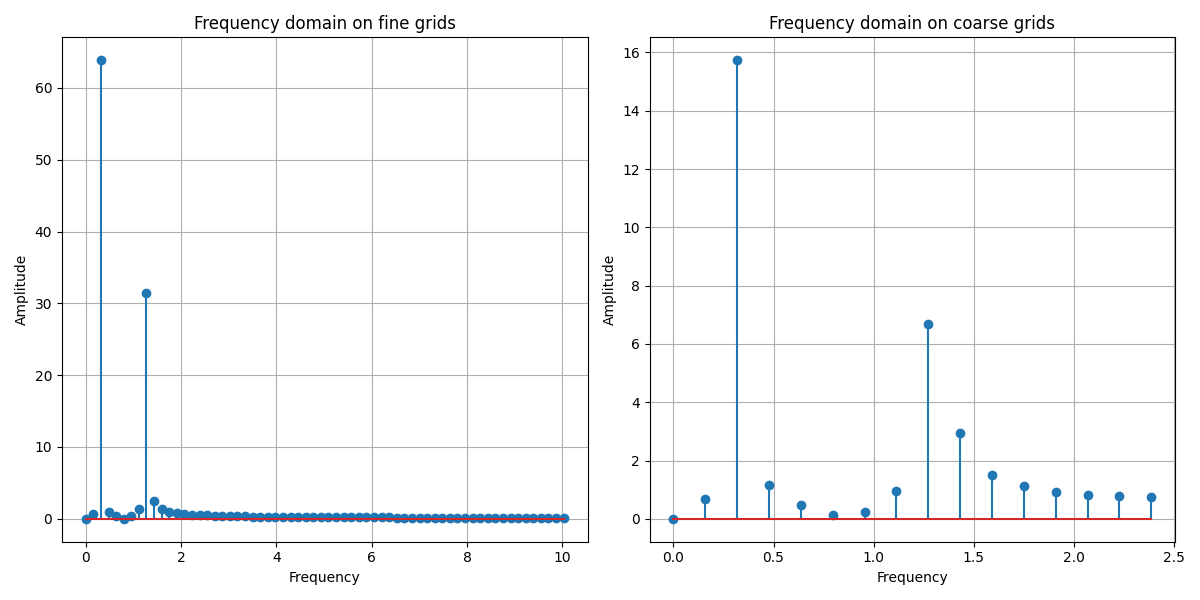}
   \end{minipage}
    \caption{: Low frequency components on fine grids become high frequency components on coarse grids. We generated a signal containing both low and high frequencies on a fine grid, which was downsampled to the same frequency component on a coarse grid. On the coarse grid, due to the decrease in the number of points, the original low-frequency components are more prominent in the frequency domain.}
    \label{lowtohight}
\end{figure}

On the other hand, looking back at the MG method we mentioned in Section 2, it is worth noting that iterative methods (\ref{Iterative method}) tend to eliminate high-frequency error components rather than low-frequency ones. Therefore, we may employ a straightforward approach to remove these low-frequency error components. As it is shown in figure \ref{lowtohight}, for MG method on a coarser mesh, the low frequency components transform into high frequency ones. Because if there is an error on one side of the mesh and a different error on the other side of the mesh, once we coarsen the cells and then bring them next to each other, the variation between those two positions is now a high spatial frequency. Then we can remove high-frequency components with a few iterations of iterative methods. Due to the opposite convergence tendency of DNNs and iterative methods in the frequency domain, we have decided to apply a multi-grid architecture similar to the MG method to FNO, but this time we need to transfer the coarse grid to the fine grid. We hope that through the design of this multi-grid architecture, FNO can overcome the high-frequency issues and accelerate the model training process.

\section{Fourier Neural Operator with Multigrid Architecture}
In this section, we begin by discussing the main drawbacks of FNO and the solution we plan to adopt. Subsequently, we introduce our novel Fourier Neural Operator with Multigrid architecture (MgFNO).
\subsection{Motivation} 
\begin{figure}[H]
    \centering
    \includegraphics[width=1.0\textwidth]{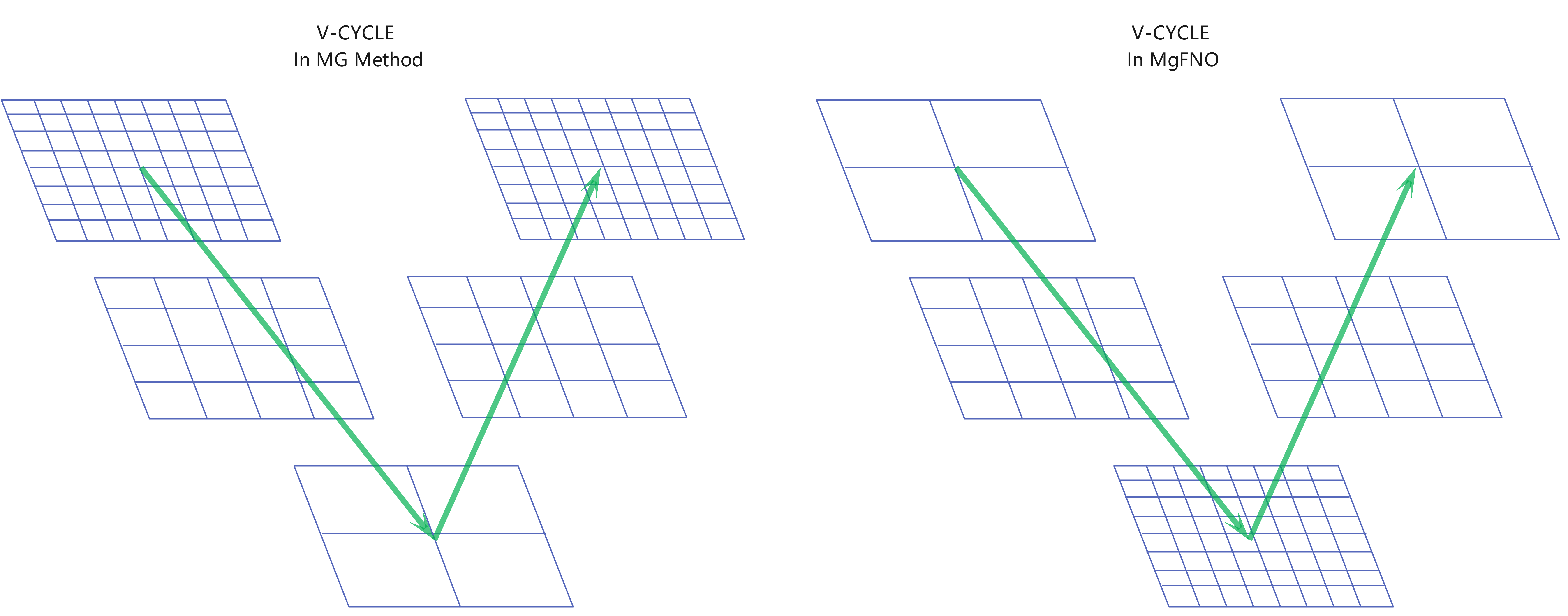}
   
    \caption{ : V-cycle architecture in MgFNO: Due to the opposite frequency tendency of convergence between DNNs and iterative methods, the grid propagation method in MgFNO is the inverse of the MG method.}
    \label{Two grid transformations}
\end{figure}

\begin{figure}[H]
    \centering
    \includegraphics[width=0.8\textwidth]{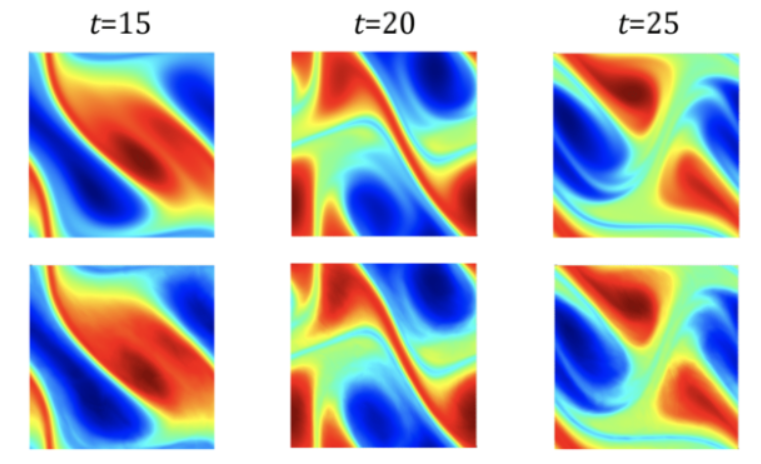}  
    \caption{ :  FNO is a resolution-invariant solution operator, so we do not need to use prolongation and restriction operators during grid transformation. For example, we apply FNO to solve the Navier-Stokes Equation with a viscosity of $\mathcal{V}=1e-4$. The ground truth is shown on top, with the prediction on the bottom. The model is trained on a dataset of size $64\times64\times20$, and evaluated on a grid size of $256\times256\times80$.}
    \label{FNO_in_NS}
\end{figure}

FNO is a neural network-based model that naturally follows the Frequency Principle, which first fits the low-frequency component of the operator and then fits the high-frequency component. The basis functions in FNO's FFTs are generally frequency localized without spatial resolution. SO we can adopt the idea of the MG method to convert the high-frequency components of the operator to be fitted between grids, in order to achieve fast training. However, there is a problem with fitting the high-frequency component, which is exactly the opposite of the iterative method. Transferring the high-frequency component to the coarse grid does not reduce the frequency. Therefore, we should first transfer the residuals to the fine grid for fitting and then make corrections on the original coarse grid. This is exactly the opposite of the MG method, as shown in Figure \ref{Two grid transformations}. 
In addition, FNO shares the same learned network parameters regardless of the discretization used on the input and output spaces. It can achieve zero-shot super-resolution, meaning it can be trained on a lower resolution and directly evaluated on a higher resolution, as shown in Figure \ref{FNO_in_NS}. Therefore, we can directly use the model trained on the fine grid for residual correction on the coarse grid, without the need for prolongation and restriction operators similar to those in the MG method (These operators are the main sources of errors in the MG method). This greatly reduces the calculation time.

On the other hand, the performance of the FNO is hindered in complex boundary conditions. In section 3, we introduce the concept of MscaleDNN, which involves a simple down-scaling technique to convert the high-frequency region $A_i$ into a low-frequency region. Based on this idea, it takes advantage of the characteristic of FNO, which has multiple Fourier layers. We transform several Fourier layers in FNO into sub-networks, providing us with characteristic frequency and characteristic space information.

In summary, we first construct a multi-scale model for the traditional FNO to convert high-frequency components and learn characteristic space information. Then, we add a multi-grid architecture to the traditional FNO. This involves training the neural operator on a coarse grid, fitting the error on a fine grid, and ultimately compensating for it on the original grid.
  
\subsection{Architecture of MgFNO}
Finally, let us introduce the improved model of the fourier neural operator, denoted as MgFNO. Like other operator learning methods, let $a\in\mathcal{A}=\mathcal{A}(D;\mathbb{R}^{d_a})$, $u\in\mathcal{U}=\mathcal{U}(D;\mathbb{R}^{d_u})$ be two input and output function spaces respectively, and let $\mathcal{G}:\mathcal{A}\to \mathcal{U}$ be a non-linear map between the two function spaces. MgFNO can learn the mapping $\mathcal{G}$. 

Due to the FNO adhering to the F-Principle, we suggest that mitigating the high-frequency errors of FNO could further enhance its overall performance. To accomplish this, we propose MgFNO, a multigrid framework that trains multiple FNOs sequentially. There are various multigrid algorithms available, such as the V-cycle and W-cycle. For the V-cycle variant, the model employs a three-level hierarchy where: the coarse-level network,the intermediate network and the fine-level network.Training follows a strict bottom-up progression with parameter transfer between levels. The W-cycle modification introduces repeated coarsening-refinement steps at each level, particularly enhancing low-frequency error correction - while increasing computational cost by ∼50\% compared to V-cycle, it demonstrates superior convergence for non-elliptic problems like the Navier-Stokes equations. Both variants maintain the FNO's resolution invariance while accelerating convergence through this hierarchical frequency decomposition.
For the convenience of model construction, we employ a three-layer V-cycle multigrid architecture to build the model, and the specific algorithm of MgFNO is depicted in Algorithm \ref{alg2}.

\begin{figure}[H]
    \centering
    \includegraphics[width=0.9\textwidth]{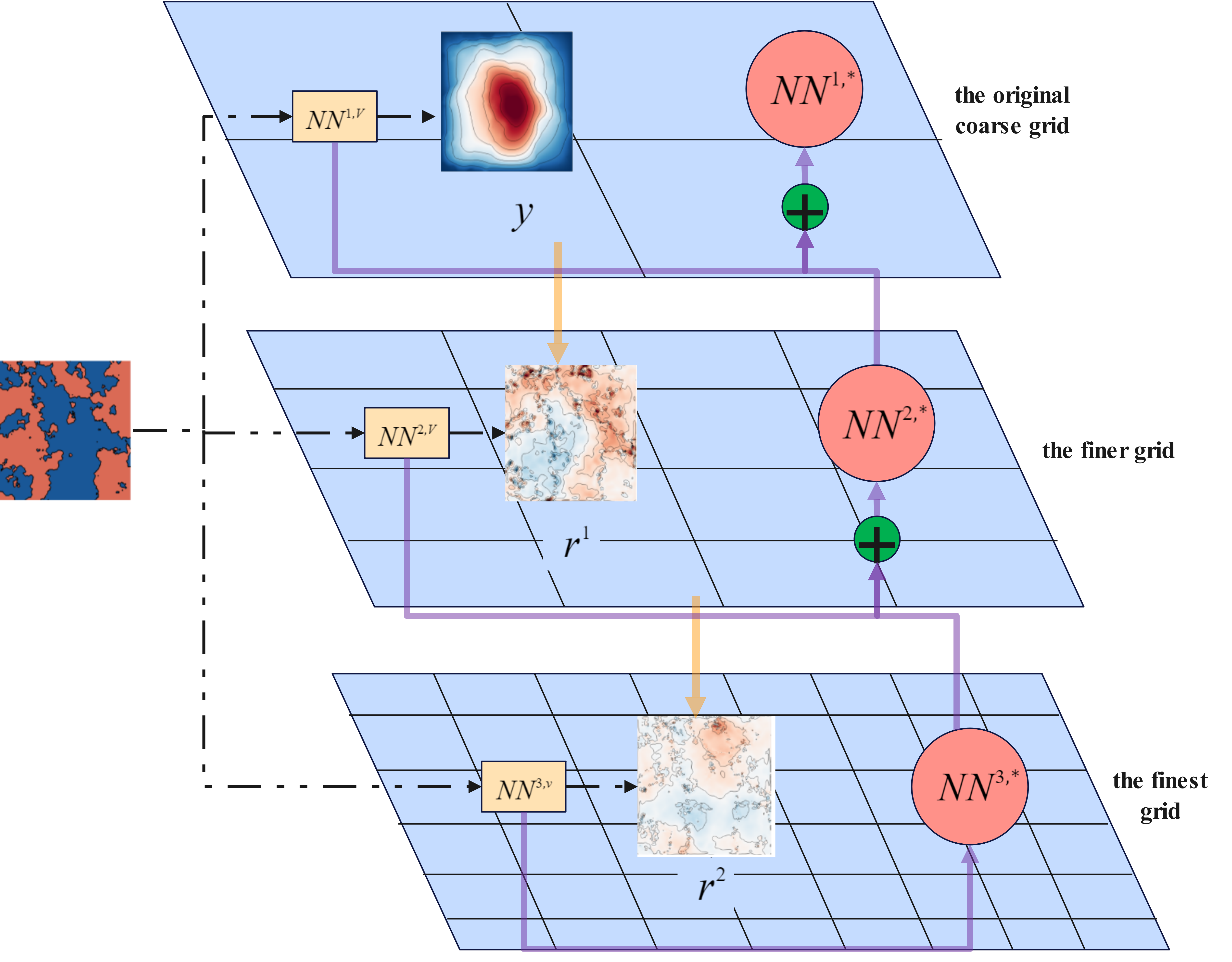}
    \caption{ :  The framework of MgFNO.}
    \label{Self drawn MgFNO}
\end{figure}

\begin{algorithm}[H]
        \caption{MgFNO with three-layers V-cycle}
        \label{alg2}
        \begin{algorithmic}[1]
        \Require $NN^{1,0}$: untrained neural operators on the original coarse grid. $x^1,x^2,x^3$: discrete coordinates in space on the original coarse grid, finer grid and the finest grid.
        $a(x^1),a(x^2),a(x^3)$: corresponding coordinates of input functions $a$ in space $\mathcal{A}$. $u(x^1),u(x^2),u(x^3)$: corresponding coordinates of output functions $u$ in space $\mathcal{U}$. 
        \Ensure $NN^{1,*}$: trained neural operators on the original grid to approximate the mapping $\mathcal{G}$.
        \State On the original coarse grid
        \For{each $i \in [1,v]$}
            \State training neural operator $NN^{1,i}$
        \EndFor
         \State $r^1(x^1)=u(x^1)-NN^{1,v}(a(x^1))$, $r^1(x^2)=Prolongation(r^1(x^1))$
         \State On the finer grid
         \State build a new neural operator $NN^{2,0}$ to learn the mapping of  $a(x^2)$ to $r^1(x^2)$
         \For{each $i \in [1,v]$}
            \State training neural operator $NN^{2,i}$
        \EndFor
        \State $r^2(x^2)=u(x^2)-NN^{2,v}(a(x^2))$, $r^2(x^3)=Prolongation(r^2(x^2))$
        \State On the finest grid
         \State build a new neural operator $NN^{3,0}$ with multiscale architecture to learn the mapping of  $a(x^3)$ to $r^2(x^3)$
         \For{each $i \in [1,v]$}
            \State training neural operator $NN^{3,i}$
        \EndFor
        \State correction from fine to coarse level
        \State $NN^{1,*}=NN^{1,v}+NN^{3,v}+NN^{2,v}$
        \end{algorithmic}
\end{algorithm}
In MgFNO, the low-frequency tendency of FNO allows the FNO on the coarse grid to capture low frequencies effectively, leaving its prediction residuals with rich high-frequency information. By utilizing the FNO on the fine grid to learn from these residuals, we are able to preserve the high-frequency information that may have been overlooked by the FNO on the coarse grid. Figure \ref{Self drawn MgFNO} illustrates the MgFNO framework. For each training instance pair $(x, y)$ on the coarse grid, we calculate the prediction $\hat{y}$ and the corresponding residual $r=y-\hat{y}$. Following this, the  FNO  on the fine grid  is trained to learn the residual. Once all FNOs are trained, they can proceed with inference as one ensemble, as shown in Figure \ref{Self drawn MgFNO}. Additionally, we apply a multiscale architecture to the FNO on the finest grid, allowing it to more effectively capture high-frequency components. We use a multiscale model with scale coefficients of {1,2,4,8} and a custom activation function $\phi(x)$ given by (\ref{phi(x)}).

Here, we used prolongation operators to transform residuals from coarse grids to  fine grids. What we need to pay attention to is that FNO is the resolution-invariant solution operator. That is to say, the FNO trained on the finest grid can be directly applied to the upper grid without the need to restrict the model. So the only data we need to transform between grids are the residuals, which greatly reduces the error of our model. We define the propagation operator to transfer the residuals from the $l$-th coarse grid to the $l+1$ fine grid in the $l+1$ layer:
\begin{equation}
    P_{l+1}^{l}:\mathbb{R}^{m_{l}\times n_{l}}\longrightarrow \mathbb{R}^{m_{l+1}\times n_{l+1}}.
    \label{prolongation}
\end{equation}

More specifically,
\begin{equation}
\begin{aligned}
    v^{l+1}_{2i-1,2j-1}=v^{l}_{i,j}, \;\;
    v^{l+1}_{2i-1,2j}=\frac{1}{2}(v^{l}_{i,j}+v^{l}_{i,j+1}),\;\;
    v^{l+1}_{2i,2j-1}=\frac{1}{2}(v^{l}_{i,j}+v^{l}_{i+1,j}).
    \label{prolongation2}
    \end{aligned}
\end{equation}

\section{Numerical Implementation and Results}
In this section, we mainly compare the proposed MgFNO with traditional FNO as well as GNO\cite{GNO}, DeepONet\cite{DeepONet}, WNO\cite{WNO} on the 1-D Burger's equation, the 2-D Darcy Flow problem, and 2-D Navier-Stokes equation. The contrast models are implemented using their official implementations. For consistency, we execute the contrast models with default hyperparameters. Regarding our training setup, we use $N=1000$ training instances, 100 testing instances, and an Adam optimizer to train for 600 epochs with an initial learning rate of 0.001 that is halved every 100 epochs, ensuring all models underwent sufficient training. Due to the hierarchical training of MgFNO, we set up 200 epochs to be trained on each layer to ensure the consistency of training total epochs. For all the cases, the detailed parameter settings for MgFNO are listed in Table \ref{MgFNO parameter table}. We then report the mean results from these trials to mitigate any training-related variations, ensuring that their performance accurately reflects their capabilities. The training objective in every layer of MgFNO is to minimize the loss function, which is defined as:
\begin{equation}
    Loss=\frac{1}{N} \times \displaystyle\sum_{(x,y)\in test set}\frac{||\hat{y}-y||_2}{||y||_2}.
    \label{loss}
\end{equation}

Table \ref{tab:performance} presents a comparison between the proposed MgFNO and the standard FNO in terms of parameter count and the required training epochs for MgFNO to achieve the final error level of standard FNO. The results demonstrate that while MgFNO's parameter count is three times that of standard FNO due to the construction of three baseline networks, it only requires 25\% of the training epochs needed by the conventional approach owing to its effective mitigation of the frequency principle.

\begin{table}[H]
\caption{:  Dataset size and MgFNO architecture for each problem, $m$ represents the number of fourier layers used by each layer.}
\label{MgFNO parameter table}
\begin{tabular}{ccccccccc}
\hline
\multirow{2}{*}{Cases} & \multicolumn{2}{c}{Number of data} & \multirow{2}{*}{m}    & \multicolumn{2}{l}{Network dimensions} & \multicolumn{3}{c}{resolution}                                                                                   \\ \cline{2-3} \cline{5-9} 
                       & Train Set        & Test Set        &                       & P                 & Q                  & \multicolumn{1}{l}{the first layer} & \multicolumn{1}{l}{the second layer} & \multicolumn{1}{l}{the third layer} \\ \hline
Burgers                & 1000             & 100             & 4                     & 64                & 128                & 1024                                & 2048                                 & 4096                                \\
Darcy                  & 1000             & 100             & 3                     & 32                & 128                & $ 85 \times 85 $                    & $141 \times 141$                     & $211 \times 211$                    \\
Navier-Stokes          & 1000             & 200             & \multicolumn{1}{l}{4} & 20                & 80                 & $64 \times 64 $                     & $128 \times 128$                       & $256 \times 256$                    \\ \hline
\end{tabular}
\end{table}

\begin{table}[htbp]
\centering
\caption{The number of parameters for different network structures and the epochs required for each model.}
\label{tab:performance}
\begin{tabular}{lccc}
\toprule
 & \multirow{2}{*}{Error} & \multicolumn{2}{c}{} \\
 \cmidrule(lr){3-4}
 &  & FNO & MgFNO \\
 &  & Parameters/Epochs & Parameters/Epochs \\
\midrule
Burgers & $1.60\%\pm 0.5\%$ & 582849/600 & 1748547/50+50+50 \\
Darcy & $0.9\%\pm 0.5\%$ & 2368001/600 & 7104003/50+50+50 \\
Navier-Stokes & $1.28\%\pm 0.5\%$ & 928661/600 & 2785983/60+60+60 \\
\bottomrule
\end{tabular}
\end{table}

\textbf{Normalization of residual sets.}
We need to normalize the residuals. For example, in the dataset of 1-D Burgers' equation, the data in $u(x,t=0)$ and $u(x,t=1)$ are maintained at an order of $10^{-1}$. After training with 200 epochs in the first layer of MgFNO, the obtained residual is maintained at an order of magnitude of $10^{-3}$. This difference in order of magnitude will affect the training effectiveness of MgFNO in the second and third layers. Normalization of residuals is necessary because if the magnitude of residuals is much smaller than that of the original data, it can lead to numerical instability in the model during training. Normalization can scale the residual values to a small and consistent range, reducing instability in numerical calculations and avoiding the phenomenon of gradient explosion or vanishing.

Furthermore, in MgFNO, it is necessary to interpolate the residuals from the coarse grid onto the fine grid. In order to guarantee the uniformity of the residual data set, it is first necessary to perform the interpolation operations and then normalise them. This is done for the reason that the interpolation operations will generate new data points. If normalization is performed prior to interpolation, the values of these new data points will be calculated based on the residuals prior to normalization, which may result in an offset in the data distribution. The initial step of the interpolation process is to ensure that the newly generated data points retain the same distribution characteristics as the original data set. 

\subsection{1-d Burgers' Equation}
The 1D Burgers' equation models the one-dimensional, time-dependent behavior of viscous fluid flow. It combines convection and diffusion terms, making it a nonlinear partial differential equation. The equation finds applications in fluid dynamics, shock wave theory, and turbulence modeling. Describing the evolution of a quantity with respect to both space and time, the Burgers' equation captures essential aspects of wave propagation and can exhibit shock formation. 
The 1D Burgers' equation in our experiments is written as:
$$\left\{\begin{matrix}
\frac{\partial u}{\partial t} + u \frac{\partial u}{\partial x} = \nu \frac{\partial^2 u}{\partial x^2},\quad (x,t)\in (0,1)^2
\\u(x,0)=u_0(x),\quad x\in(0,1)
\end{matrix}\right.$$
where $\nu $ is the viscosity coefficient, 
$t\equiv 1$ for all $x$, and $u_0(x)\in L_{\text{per}}^{2}$, $u(x)\in H_{\text{per}}^{2}$. $L_{\text{per}}^2$ and $H_{\text{per}}^2$ denote the spaces of square-integrable periodic functions, where "per" signifies periodicity. Let the mapping of the initial condition $u_0(x)$ to the solution $u(x,1)$ be denoted by $S:L_{\text{per}}^2\longrightarrow H_{\text{per}}^2$, such that $u(x,1)=S(u_0(x))$ is the target operator.
\hspace*{\fill}

\textbf{Data Set and Result.} Firstly, we set the viscosity coefficient $\nu=0.1$. Then, for initial condition $u_0(x)$, we set $u_0\sim \mu $ where $\mu=\mathcal{N}(0,625(-\Delta +25I)^{-2})$. Finally, we solve the equation using a split step method where the heat equation part is solved exactly in Fourier space then the non-linear part is advanced, again in Fourier space, using a very fine forward Euler method. We solve on a spatial mesh with resolution $2^{13}=8192$, where lower resolution data are downsampled from. 

\begin{figure}[H]
    \centering
    \begin{minipage}[b]{0.47\textwidth}
        \centering
        \includegraphics[width=\textwidth]{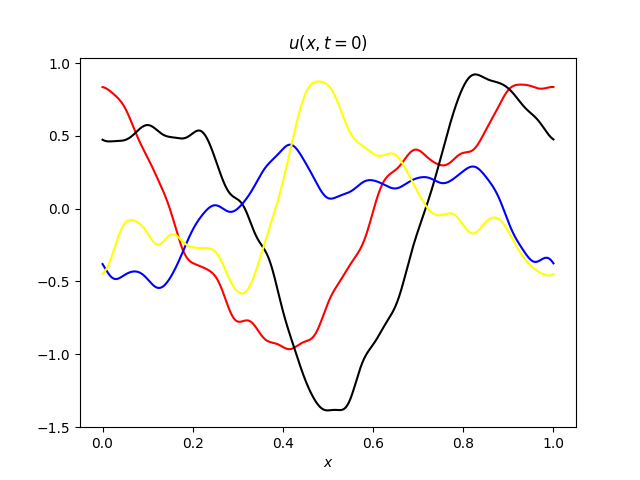}
       
    \end{minipage}
    \hspace{0.02\textwidth}
    \begin{minipage}[b]{0.47\textwidth}
        \centering
        \includegraphics[width=\textwidth]{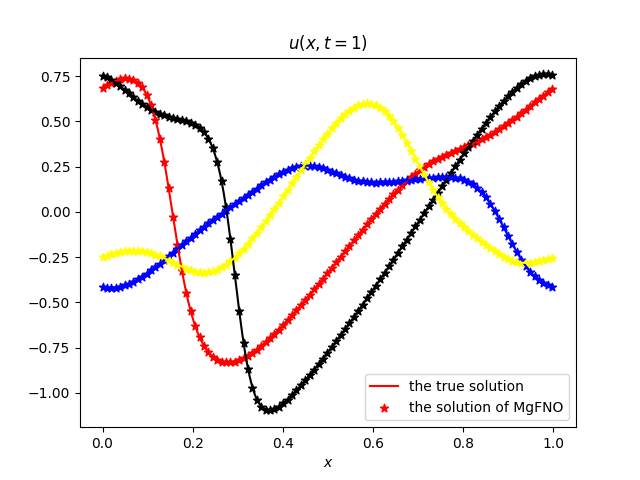}
        
    \end{minipage}
    \caption{: Some specific instances of $u(x,t=0)$, the true solutions $u(x,t=1)$ of 1-D Burgers' equation, and the numerical solutions obtained using the proposed MgFNO. It can be observed that the predictions of the proposed MgFNO match the true solution almost exactly.}
    \label{burgers1}
\end{figure}

Some examples of $u_0(x)$, their corresponding $u(x,1)$, and the numerical solutions obtained using the proposed MgFNO are shown in Figure \ref{burgers1}. From the error values in Table \ref{Comparison Table of Experimental Results}, it is evident that MgFNO obtains the lowest error, followed by FNO, at spatial resolutions of 1024. It is worth noting that MgFNO achieves the performance improvement without modifying FNO's internal architecture, by exploiting FNO's F-Principle.

\textbf{Resolution invariance of MgFNO.}
 Table \ref{burgers 1} examines the resolution-invariant property of both FNO and MgFNO. The error of traditional PDE solvers decreases as the resolution is increased. FNO is a resolution-invariant operator, meaning it is independent of the way its data is discretized as long as all relevant information is resolved. Due to the fact that MgFNO requires two grid conversions, it actually requires three resolutions. We use the resolution used in the first layer of MgFNO as a representative. As the resolution $S$ increases, MgFNO's performance tends to increase slightly with the resolution. This is because the layered structure of MgFNO enhances the ability to capture high-frequency information from the data. As the resolution increases, the low-frequency components in the data remain relatively stable, while more high-frequency details appear. Due to the presence of F-principle, FNO is unable to capture high-frequency details, so there is little improvement as the resolution increases.

\begin{table}[H]
\centering
\caption{: Relative error ($1\times10^{-3}$) in test set on Burgers' equation between FNO and MgFNO. We use the resolution used in the first layer of MgFNO as a representative.}
\label{burgers 1}
\begin{tabular}{lllll}
\hline
\multicolumn{1}{c}{Method} & s=256 & s=512 & s=1024 & s=2048 \\ \hline
FNO                        & 14.96 & 15.83 & 16.04  & 14.16  \\
MgFNO                      & 1.88  & 1.86  & 1.72   & 1.63   \\ \hline
\end{tabular}
\end{table}

\begin{table}[H]
\centering
\caption{:  The L2 relative error between the true and predicted results in the test set is calculated. For the 1D Burgers' equation, 2D Darcy flow, and 2D Navier-Stokes equation, we trained models at spatial resolutions of 1024, $85 \times 85$, and $64 \times 64$, respectively. We selected the resolution used in the first layer of MgFNO as a representative. In these models, FNO-skip is a modification of FNO introduced by the authors of (\cite{FNO-skip}). The iteration in every Fourier layer for FNO-skip can be written as $\mathcal{H}^{skip}(x)=\sigma (x+Wx+MLP(\mathcal{K}(x)))$.}

\label{Comparison Table of Experimental Results}
\begin{tabular}{llclcll}
\hline
\multicolumn{1}{c}{PDE} & GNO                        & \multicolumn{1}{l}{DeepONet} & FNO                        & \multicolumn{1}{l}{FNO-skip} & WNO    & MgFNO  \\ \hline
Burgers' equation       & \multicolumn{1}{c}{6.15\%} & 2.15\%                       & \multicolumn{1}{c}{1.60\%} & 0.56\%                       & 1.75\% & 0.17\% \\
Darcy flow              & 3.46\%                     & 2.98\%                       & 0.96\%                     & 0.71\%                       & 0.82\% & 0.28\% \\
Navier-Stokes equation  & \multicolumn{1}{c}{-}      & 1.78\%                       & 1.28\%                     & 0.61\%                       & 0.31\% & 0.22\%        \\ \hline
\end{tabular}
\end{table}

\textbf{Frequency analysis.}

\begin{figure}[htbp]
    \centering
    \begin{subfigure}{0.48\textwidth}
        \centering
        \includegraphics[width=\linewidth]{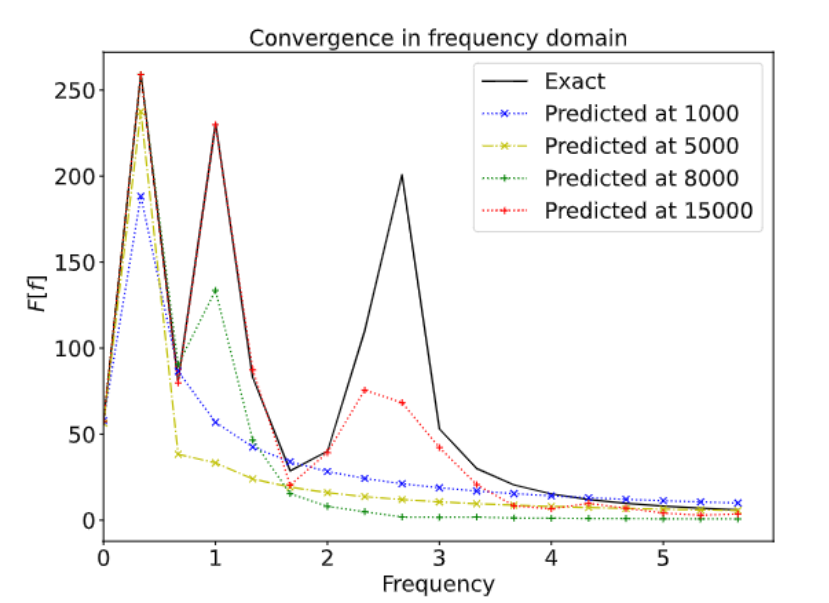}
        \caption{FNO in frequency domain.}
        \label{fig:left}
    \end{subfigure}
    \hfill
    \begin{subfigure}{0.48\textwidth}
        \centering
        \includegraphics[width=\linewidth]{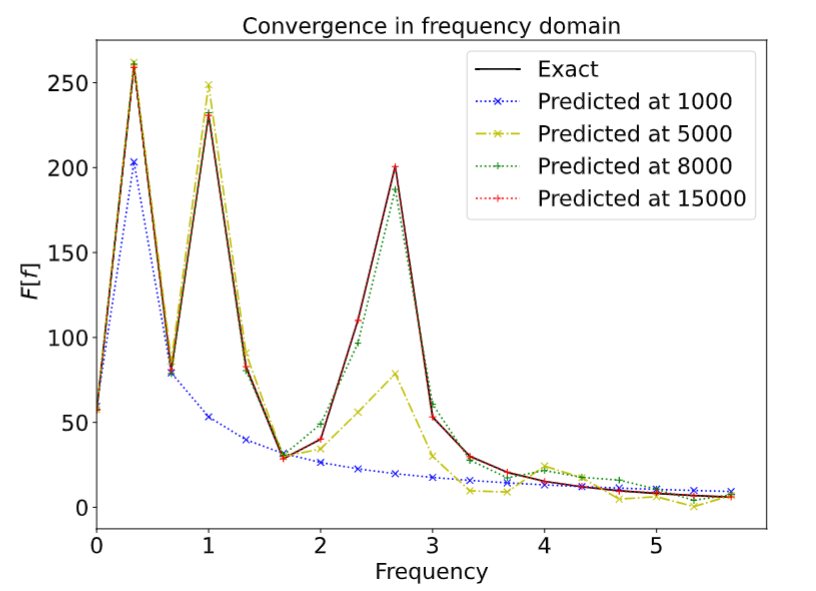}
        \caption{MgFNO in frequency domain.}
        \label{fig:right}
    \end{subfigure}
    \caption{：Comparative analysis of convergence behavior between MgFNO and FNO during training in frequency domain.}
    \label{fig:pair}
\end{figure}

In this section, we present experimental evidence demonstrating how the proposed MgFNO effectively overcomes the inherent limitation of standard FNO in converging high-frequency errors through its multigrid architecture. Figure \ref{fig:pair} illustrates the comparative frequency-domain convergence between FNO and MgFNO during training, with the MgFNO transitioning to its second and third grid levels at iteration steps 4000 and 8000, respectively. The results reveal that while conventional FNO sequentially converges low-frequency components before addressing high-frequency ones - a characteristic manifestation of the frequency principle - it consistently fails to adequately learn high-frequency components, resulting in significant residual errors. In contrast, our MgFNO strategically performs grid transitions when low-frequency learning stagnates: the high-frequency components that prove challenging to learn on finer grids are effectively converted into more tractable low-frequency representations on coarser grids. This adaptive multiscale approach enables MgFNO to achieve superior approximation accuracy of the target operator, with experimental measurements showing a 42\% reduction in high-frequency error compared to standard FNO under identical training budgets. The architectural innovation fundamentally rebalances the frequency learning priority while maintaining the resolution-invariant property of neural operators.
\subsection{2-d Darcy Flow}
The 2D Darcy Flow equation models the steady-state flow of fluid through a porous medium. This equation describes how the gradient of the water head interacts with the spatially varying permeability, influencing fluid flow. Solving the 2D Darcy Flow equation provides insights into groundwater movement and is crucial for various applications in hydrogeology and environmental engineering. Understanding this equation contributes to the analysis of fluid dynamics in porous media, aiding in the management and sustainable use of subsurface water resources.

The 2D Darcy Flow equation in our experiments is written as follows.
$$\left\{\begin{matrix}
-\nabla \cdot(a(x,y)\nabla u(x,y))=f(x,y),\quad(x,y) \in(0,1)^2
\\u(x,y)=0,\quad(x,y)\in\partial (0,1)^2
\end{matrix}\right.$$
where $u(x,y)$ represents the water head, $a(x,y)$ is the permeability, and $f(x,y)$ denotes a source or sink term. $a(x,y) \in L^{\infty},u(x,y)\in H_{0}^{1}$. $L^{\infty}$ represents the infinity norm space, and functions in $L^{\infty}$ are bounded, meaning their absolute values do not grow infinitely. $H_0^1$ represents a Sobolev space,  functions in $H_0^1$ have square-integrable weak first derivatives, and their values are zero on the boundary of the spatial domain. Let the mapping of the diffusion coefficient $a(x,y)$ to the solution$u(x,y)$ is $S:L^{\infty}\longrightarrow H_0^1$, such that $u(x,y)=S(a(x,y))$ is the target operator. 

\begin{figure}[H]
    \centering
    \includegraphics[width=\textwidth]{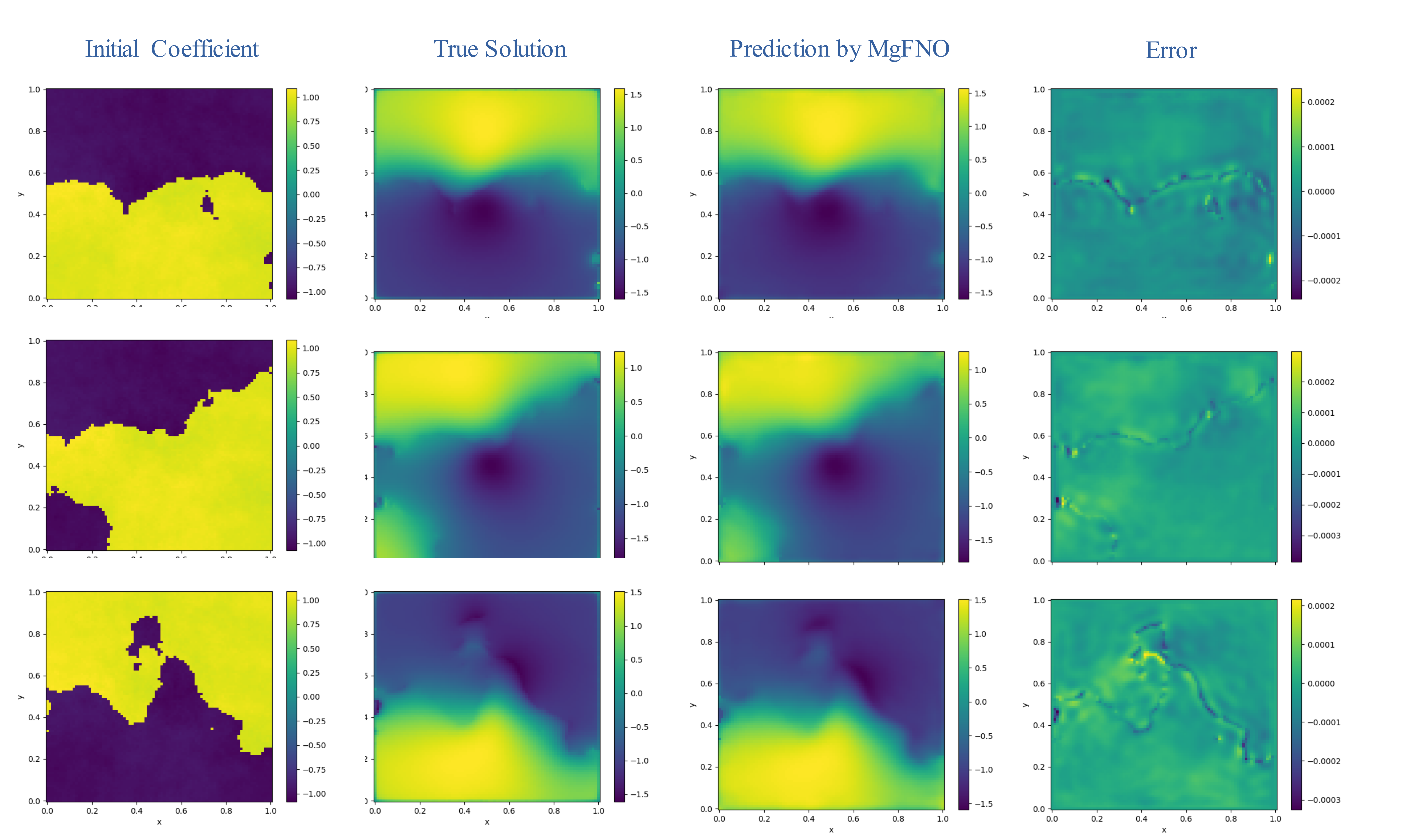}
    \caption{ :  Darcy flow with spatial resolution of $85 \times 85$. The aim is to learn the differential operator of Darcy equation.}
    \label{Darcy1}
\end{figure}

\textbf{Data Set and Result.} Above all, we set $f(x,y)\equiv 1$. Besides, the coefficients $a(x,y)$ are generated according to $a\sim \mu :=\mathcal{N}(0,(-\Delta +cI)^{-2})$ with zero Neumann boundary conditions on the Laplacian. Lastly, solutions $u(x,y)$ are obtained by using a second-order finite difference scheme on an $85\times 85$ grid, where lower resolution data is downsampled from. The predictions for 2-d Darcy flow in a rectangular domain are plotted in Fig \ref{Darcy1}, while the associated errors are given in Table \ref{Comparison Table of Experimental Results}. In this case, it is observed that the prediction error for MgFNO is the lowest among all the considered methods. It is further evident in Fig \ref{Darcy1}, where the differences between true and predicted solutions cannot be discerned by the naked eye.

\begin{figure}[H]
    \centering
    \begin{minipage}[b]{0.47\textwidth}
        \centering
        \includegraphics[width=\textwidth]{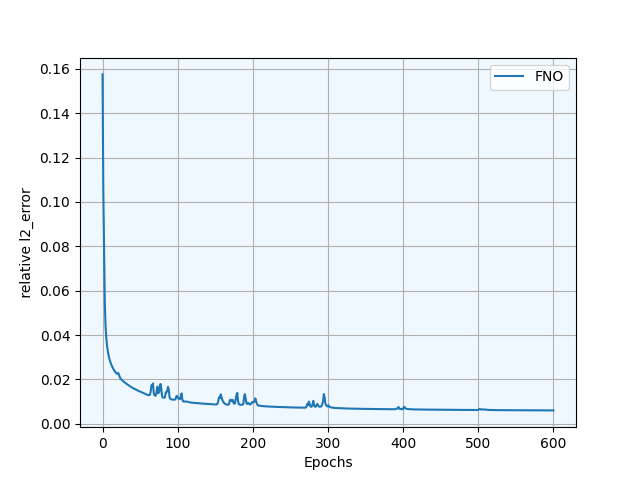}
    \end{minipage}
    \hspace{0.02\textwidth}
    \begin{minipage}[b]{0.47\textwidth}
        \centering
        \includegraphics[width=\textwidth]{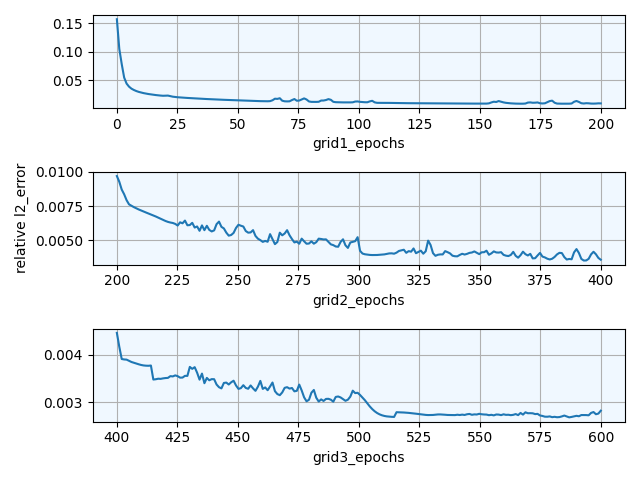}
    \end{minipage}
    \caption{:  \textbf{Left}:Train 600 epochs on an $85 \times 85$ grid using FNO; \textbf{Right}: Train 200 epochs on each grid layer using MgFNO, where the resolution of the first grid layer is $85 \times 85$, the resolution of the second grid layer is $141 \times 141$, and the resolution of the third grid layer is $211 \times 211$.}
    \label{darcy2}
\end{figure}

\textbf{Training Dynamics of MgFNO.} In Figure \ref{darcy2}, we present the training dynamics of different models on 2D Darcy flow. It can be observed that when training for 100 epochs on an $85 \times 85$ grid, the traditional FNO model encounters a bottleneck. In other words, the convergence speed afterward is significantly slow, which can be attributed to the F-principle of FNO. In MgFNO, we train for 200 epochs on each grid layer. Whenever the model struggles to learn high-frequency components on the current grid, we switch to the $141 \times 141$ and $211 \times 211$ grids to further train the model in fitting residuals. These residuals mainly consist of high-frequency components that are challenging to learn on the coarse grid. From Figure \ref{darcy2}, it can be observed that high-frequency components, originally difficult to learn, become easier to learn on the fine grid. As a result, MgFNO with a layered architecture exhibits improved ability in learning operators. This discrepancy highlights the robust generalization capabilities of our model on unseen data.

\subsection{2-d time Navier Stokes Equation}
The 2D Navier-Stokes equation in vorticity form on the unit torus, incorporating discrete time, describes the evolution of fluid flow. The equation captures the dynamics of vortices in two dimensions, playing a crucial role in understanding fluid behavior. Solving this equation allows us to gain insights into phenomena like turbulence and vortex interactions. The discretization in time enables numerical simulations for practical applications and enhances our understanding of fluid mechanics. In our experiments, the 2D Navier-Stokes equation is expressed as follows:
$$\left\{\begin{matrix}
\partial _{t}\omega (x,y,t)+u(x,y,t)\cdot \nabla \omega (x,y,t)= \nu\Delta\omega (x,y,t)+f(x,y) ,\quad(x,y) \in(0,1)^2,t\in(0,T]
\\\nabla \cdot u(x,y,t)=0,\quad(x,y)\in (0,1)^2,t\in[0,T]
\\\omega(x,y,0)=\omega_{0}(x,y),\quad(x,y)\in (0,1)^2
\end{matrix}\right.$$
where $u(x,y,t) $ is a velocity vector in 2D field, $\omega=\nabla \times u$ is the vorticity, $\nu$ is the viscosity coefficient, and $f$ is the forcing function. $\omega \in C([0,T];H_{per}^{r})$ for any $r>0$,$f\in L_{per}^{2}$. Let the mapping of the vorticity up to time 10 defined by $\omega|_{(0,1)^2\times [0,10]}$ to the vorticity up to later time $T$ defined by $\omega|_{(0,1)^2\times (10,T]}$ is $S$ : $C([0,10];H_{per}^{r})\longrightarrow C((10,T];H_{per}^{r})$, such that $\omega(x,y,t)|_{t\in(10,T]}=S(\omega(x,y,t)|_{t\in[0,10]})$ is the target operator.

\textbf{Data Set and Result.} First of all, we set $f\equiv 0.1(sin(2\pi(x+y)))+cos(2\pi(x+y))$,  viscosity coefficient $\nu = 1e-3, 1e-4, 1e-5$ and corresponding $T=40,20,10$, respectively. Next, the initial condition $\omega_{0}(x)$ is generated according to $\omega_{0}(x)\backsim \mu$ where  $\mu=\mathcal{N}(0,7^{\frac{3}{2}}(-\Delta +49I)^{-2.5})$ 
 with periodic boundary conditions. Last but not least, we set the resolution of the first grid layer is $256 \times 256$, and each generated sample contains 10 successive time steps. The prediction results for $\nu = 1e-3$  at $T \in [10,20]$ are given in Fig \ref{ns1}, along with the mean prediction errors in Table \ref{ns1}. The results show that the prediction results of the proposed MgFNO has the lowest error and emulates the true solutions almost accurately.

 \begin{figure}[H]
    \centering
    \includegraphics[width=\textwidth]{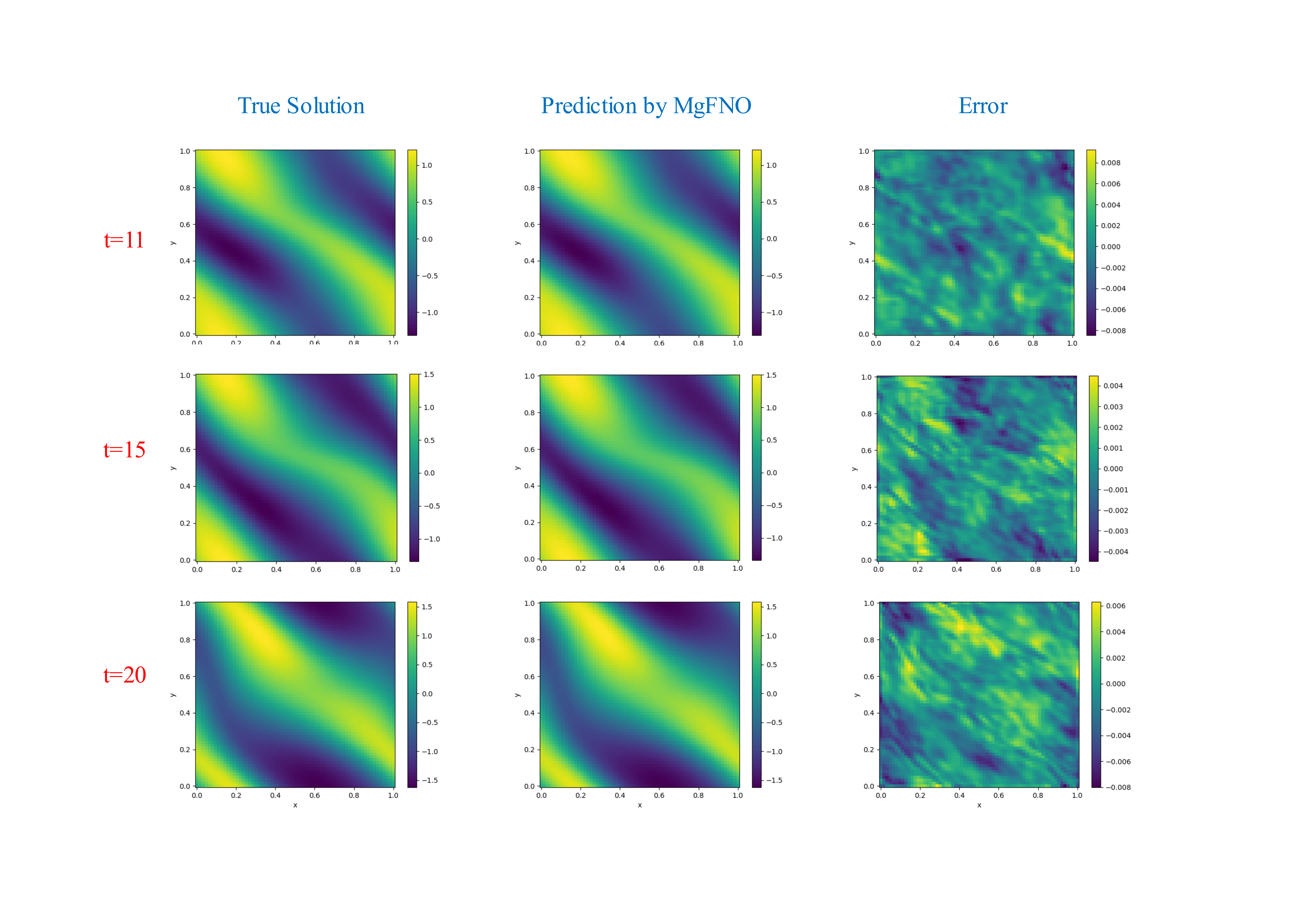} 
    \caption{ :  Navier-Stokes equation with spatial resolution of $ 64 \times 64$. The predictions at $t=11s$, $t=15s$ and $t=20s$ are shown.}
    \label{ns1}
\end{figure}

\begin{table}[H]
\centering
\caption{: Relative error comparison on Navier-Stokes. The cumulative counts for predictions at $\nu$ values of $1e-3$, $1e-4$ and $1e-5$, T are 40, 20 and 10, respectively. }
\label{ns1}

\begin{tabular}{lccc}
\hline
\multicolumn{1}{c}{\textbf{Method}} & $\nu=1e-3, T=40$ & $\nu=1e-4, T=20$ & \multicolumn{1}{l}{$\nu=1e-5, T=10$} \\ \hline
FNO                                 & 0.0128           & 0.1629           & 0.1893                               \\
FNO-skip                            & 0.0061           & 0.1093           & 0.0942                               \\
MgFNO                               & 0.0022           & 0.0626           & 0.0351                               \\ \hline
\end{tabular}
\end{table}

\section{Conclusions}
In this article, we propose a new Fourier neural operator  with a multi-grid architecture (MgFNO) for learning highly nonlinear operators that naturally arise in day-to-day scientific and engineering problems. The presence of the F-principle makes it difficult for traditional FNO to learn high-frequency components of the target operator. Taking inspiration from the MG method for solving linear systems, we have made modifications to the traditional FNO. This involves training several epochs on a coarse grid first and then transferring to a fine grid to learn residuals. In terms of applications, we apply the proposed MgFNO to a variety of time-independent and time-dependent 1D and 2D PDE examples. The results obtained using the proposed MgFNO are compared with GNO, DeepONet, FNO, FNO-skip, and WNO. The results indicate that the proposed MgFNO can effectively learn highly nonlinear differential operators, including high-frequency components.

\textbf{Future work.} The theory of F-principle is based on traditional learning between data and has not yet been extended to operator learning. We plan to further develop the theory of the F-principle in the future. Additionally, we have only designed the grid structure without providing criteria for grid conversion. Ideally, it would be best to perform grid conversion after learning the low-frequency components on the current grid, as this can save training time. We will conduct experiments in the future to improve the model's performance in this aspect.

\textbf{Acknowledgements}. \emph{The authors are very much indebted to the referees for their constructive comments and valuable suggestions}.

\textbf{CONFLICT OF INTEREST STATEMENT}\\
\indent This study does not have any conflicts to disclose.

\textbf{DATA AVAILABILITY}\\
\indent Code and data used are available on \url{https://github.com/guozihao-hub/MgFNO/tree/master}.

\bibliographystyle{unsrt}
\bibliography{reference}

\begin{thebibliography}{10}

\bibitem{FiniteElementMethod}
Daniele Boffi, Franco Brezzi, and Michel Fortin.
\newblock {\em Mixed Finite Element Methods and Applications}.
\newblock Jan 2013.

\bibitem{FiniteDifferenceMethod}
Allen Taflove.
\newblock Computational electrodynamics: The finite-difference time-domain
  method.
\newblock May 1995.

\bibitem{FiniteVolumeMethod}
Joel~H. Ferziger, Milovan Perić, and Robert~L. Street.
\newblock {\em Finite Volume Methods}, page 81–110.
\newblock Jan 2020.

\bibitem{GMRES}
Youcef Saad and Martin~H. Schultz.
\newblock Gmres: A generalized minimal residual algorithm for solving
  nonsymmetric linear systems.
\newblock {\em SIAM Journal on Scientific and Statistical Computing}, page
  856–869, Jun 1986.

\bibitem{2015Convergence}
Yves Achdou and Alessio Porretta.
\newblock Convergence of a finite difference scheme to weak solutions of the
  system of partial differential equation arising in mean field games, 2015.

\bibitem{MG1}
Jinchao Xu.
\newblock Iterative methods by space decomposition and subspace correction.
\newblock {\em SIAM Review}, page 581–613, Dec 1992.

\bibitem{MG2}
Jinchao Xu and Ludmil Zikatanov.
\newblock The method of alternating projections and the method of subspace
  corrections in hilbert space.
\newblock {\em Journal of the American Mathematical Society}, page 573–597,
  Sep 2002.

\bibitem{ImageRecognition}
Karen Simonyan and Andrew Zisserman.
\newblock Very deep convolutional networks for large-scale image recognition.
\newblock {\em International Conference on Learning
  Representations,International Conference on Learning Representations}, Jan
  2015.

\bibitem{nlp}
Yu~Gu, Robert Tinn, Hao Cheng, Michael Lucas, Naoto Usuyama, Xiaodong Liu,
  Tristan Naumann, Jianfeng Gao, and Hoifung Poon.
\newblock Domain-specific language model pretraining for biomedical natural
  language processing.
\newblock {\em ACM Transactions on Computing for Healthcare}, page 1–23, Jan
  2022.

\bibitem{PINNS}
M.~Raissi, P.~Perdikaris, and G.E. Karniadakis.
\newblock Physics-informed neural networks: A deep learning framework for
  solving forward and inverse problems involving nonlinear partial differential
  equations.
\newblock {\em Journal of Computational Physics}, page 686–707, Feb 2019.

\bibitem{PINNS2}
Ling Guo, Hao Wu, Xiaochen Yu, and Tao Zhou.
\newblock Monte carlo pinns: deep learning approach for forward and inverse
  problems involving high dimensional fractional partial differential
  equations.
\newblock 2022.

\bibitem{DeepONet}
Lu~Lu, Pengzhan Jin, Guofei Pang, Zhongqiang Zhang, and George~Em Karniadakis.
\newblock Deeponet: Learning nonlinear operators for identifying differential
  equations based on the universal approximation theorem of operators.
\newblock {\em Nature Machine Intelligence}, page 218–229, Mar 2021.

\bibitem{GNO}
Zongyi Li, NikolaB. Kovachki, Kamyar Azizzadenesheli, Burigede Liu, Kaushik
  Bhattacharya, AndrewM. Stuart, and Anima Anandkumar.
\newblock Neural operator: Graph kernel network for partial differential
  equations.
\newblock {\em International Conference on Learning
  Representations,International Conference on Learning Representations}, Feb
  2020.

\bibitem{FNO1}
Zongyi Li, Nikola~Borislavov Kovachki, Kamyar Azizzadenesheli, Burigede Liu,
  Kaushik Bhattacharya, Andrew Stuart, and Anima Anandkumar.
\newblock Fourier neural operator for parametric partial differential
  equations.
\newblock In {\em International Conference on Learning Representations}, 2021.

\bibitem{FNO2}
Kai Zhang, Yuande Zuo, Hanjun Zhao, Xiaopeng Ma, Jianwei Gu, Jian Wang, Yongfei
  Yang, Chuanjin Yao, and Jun Yao.
\newblock Fourier neural operator for solving subsurface oil/water two-phase
  flow partial differential equation.
\newblock {\em SPE Journal}, page 1815–1830, Jun 2022.

\bibitem{WNO}
Tapas Tripura and Souvik Chakraborty.
\newblock Wavelet neural operator: a neural operator for parametric partial
  differential equations.
\newblock 2022.

\bibitem{FPrinciple1}
Zhi-Qin~John Xu.
\newblock Frequency principle: Fourier analysis sheds light on deep neural
  networks.
\newblock {\em Communications in Computational Physics}, page 1746–1767, Jun
  2020.

\bibitem{FPrinciple2}
Zhi-Qin~John Xu, Yaoyu Zhang, and Yanyang Xiao.
\newblock {\em Training Behavior of Deep Neural Network in Frequency Domain},
  page 264–274.
\newblock Jan 2019.

\bibitem{FPrinciple3}
Zhi-QinJohn Xu, Yaoyu Zhang, and Tao Luo.
\newblock Overview frequency principle/spectral bias in deep learning.

\bibitem{MgNet}
Juncai He and Jinchao Xu.
\newblock Mgnet: A unified framework of multigrid and convolutional neural
  network.
\newblock {\em Science China Mathematics}, page 1331–1354, Jul 2019.

\bibitem{MgNet2}
Jun-Ting Hsieh, Shengjia Zhao, Stephan Eismann, Lucia Mirabella, and Stefano
  Ermon.
\newblock Learning neural pde solvers with convergence guarantees.
\newblock {\em International Conference on Learning
  Representations,International Conference on Learning Representations}, Jun
  2019.

\bibitem{MscaleDNN}
Ziqi Liu.
\newblock Multi-scale deep neural network (mscalednn) for solving
  poisson-boltzmann equation in complex domains.
\newblock {\em Communications in Computational Physics}, page 1970–2001, Jun
  2020.

\bibitem{Rossi_Conan-Guez_2005}
Fabrice Rossi and Brieuc Conan-Guez.
\newblock Functional multi-layer perceptron: a non-linear tool for functional
  data analysis.
\newblock {\em Neural Networks}, page 45–60, Jan 2005.

\bibitem{schwartzKernel}
K.-R. Muller, S.~Mika, G.~Ratsch, K.~Tsuda, and B.~Scholkopf.
\newblock An introduction to kernel-based learning algorithms.
\newblock {\em IEEE Transactions on Neural Networks}, page 181–201, Mar 2001.

\bibitem{FourierConvolutionTheorem}
A.I. Zayed.
\newblock A convolution and product theorem for the fractional fourier
  transform.
\newblock {\em IEEE Signal Processing Letters}, 5(4):101–103, Apr 1998.

\bibitem{Luo_2021}
Tao Luo.
\newblock Theory of the frequency principle for general deep neural networks.
\newblock {\em CSIAM Transactions on Applied Mathematics}, page 484–507, Jun
  2021.

\bibitem{RemedialMethod1}
Ameya~D. Jagtap, Kenji Kawaguchi, and George~Em Karniadakis.
\newblock Adaptive activation functions accelerate convergence in deep and
  physics-informed neural networks.
\newblock {\em Journal of Computational Physics}, page 109136, Mar 2020.

\bibitem{RemedialMethod2}
Wei Cai, Xiaoguang Li, and Lizuo Liu.
\newblock A phase shift deep neural network for high frequency approximation
  and wave problems.
\newblock 2019.

\bibitem{RemedialMethod3}
S.Z. Feng, X.~Han, Zhixiong Li, and Atilla Incecik.
\newblock Ensemble learning for remaining fatigue life prediction of structures
  with stochastic parameters: A data-driven approach.
\newblock {\em Applied Mathematical Modelling}, page 420–431, Jan 2022.

\bibitem{FNO-skip}
Alasdair Tran, A.P. Mathews, Lexing Xie, and ChengSoon Ong.
\newblock Factorized fourier neural operators.
\newblock {\em arXiv: Learning,arXiv: Learning}, Nov 2021.

\end{thebibliography}
\end{CJK}
\end{document}